\documentclass[a4paper,10pt]{article}
\usepackage{latexsym}
\usepackage{amssymb}
\usepackage{amsmath}
\usepackage{color}
\usepackage{graphics}
\newlength{\cqfd}
\newcommand{\vep}{\varepsilon}

\renewcommand{\ni}{\noindent}

\newtheorem{Theorem}{Theorem}
\newtheorem{Lemma}[Theorem]{Lemma}
\newtheorem{prop}[Theorem]{Proposition}
\newtheorem{df}{Definition}
\newtheorem{rmk}{Remark}
\newcommand{\pf}{{\noindent\it Proof.~}}

\newcommand{\disp}{\displaystyle}
\newcommand{\un}{\mathbf{1}}
\newcommand{\rs}{\rho^*}
\newcommand{\Div}{{\rm div}}
\newcommand{\na}{\nabla}
\newcommand{\re}{\rho_\vep}

\newcommand{\ue}{\vu_\vep}
\newcommand{\rd}{\rho_\delta}
\newcommand{\ud}{\vu_\delta}

\newcommand{\vp}{\varphi}

\newcommand{\pt}{\partial_{t}}
\newcommand{\px}{\partial_{x}}

\newcommand{\Dt}{\frac{ d}{dt}}
\newcommand{\dx}{{\rm d} {x}}
\newcommand{\dt}{{\rm d} t }
\newcommand{\intO}[1]{\int_{\Omega} #1 \ \dx}
\newcommand{\intOB}[1]{\int_{\Omega} \left( #1 \right) \ \dx}
\newcommand{\intT}[1]{\int_0^T #1 \ \dt}
\newcommand{\dxdt}{\dx \ \dt}
\newcommand{\intTO}[1]{ \int_0^T\!\!\!\! \int_{\Omega}  #1 \ \dxdt}

\newcommand{\eq}[1]{\begin{equation}
\begin{split}
#1
\end{split}
\end{equation}}
\newcommand{\eqh}[1]{\begin{equation*}
\begin{split}
#1
\end{split}
\end{equation*}}
\newcommand{\vc}[1]{{\bf #1}}
\newcommand{\vu}{\vc{u}}
\newcommand{\Ov}[1]{\overline{#1}}

\newcommand{\lr}[1]{\left( #1 \right)}
\newcommand{\bom}{\partial\Omega}
\newcommand{\lap}{\Delta}
\newcommand{\Pd}{\pi_{\kappa,\delta}\lr{\dfrac{\rd}{\rs}}}
\newcommand{\Pn}{\pi_{\vep}\lr{\frac{\re}{\rs}}}

\title{\bf Free/Congested Two-Phase Model from Weak Solutions to Multi-Dimensional Compressible Navier-Stokes Equations}

\begin{document}
 \maketitle 
 \begin{center}
{\large 
Charlotte Perrin \footnote{Universit\'e de Savoie, Laboratoire de Math\'ematiques
 UMR CNRS 5127, Campus Scientifique, 73376 Le Bourget du Lac, France; charlotte.perrin@univ-savoie.fr}\quad and\quad 
Ewelina Zatorska  \footnote{CMAP UMR 7641 Ecole Polytechnique CNRS, Route de Saclay, 91128 Palaiseau Cedex France; ewelina.zatorska@cmap.polytechnique.fr}
\footnote{Institute of Applied Mathematics and Mechanics,  Univeristy of Warsaw,  ul. Banacha 2, 02-097 Warszawa Poland; e.zatorska@mimuw.edu.pl} }
\end{center}
\vspace{0.2cm}
\vspace {0.5cm}
\begin{center}
{\bf Abstract}
\end{center}
We approximate a  two--phase model by the compressible Navier-Stokes equations with a singular pressure term.  Up to a subsequence, these solutions are shown to converge to a global weak solution of the compressible system with the congestion constraint studied for instance by {\sc P.--L. Lions} and {\sc N. Masmoudi} [{\it Annales I.H.P.},  1999]. The paper  is an extension of the previous result obtained in one-dimensional setting by {\sc D. Bresch} {\it et al.} [{\it C. R. Acad. Sciences Paris}, 2014] to the multi-dimensional case with heterogeneous barrier for the density. 

\section{Introduction}

Macroscopic models of moving crowd identify the swarm through some density that is transported by a velocity vector field, see for instance a review paper by {\sc B. Maury} \cite{Ma}. For example, to describe the traffic jams, one may use the one-dimensional fluid model describing two-phase flow 
\eq{
\left.
\begin{array}{r}
\pt\alpha+\px(\alpha u)=0\\
\pt(\alpha u)+\px(\alpha u^2+\pi)=0
\end{array}
\right\}
\quad\text{in } (0,T)\times \mathbb{R}\label{s0}}
with the following restrictions
\eq{0\leq\alpha\leq 1,\qquad (1-\alpha)\pi=0\label{cs0}.}
Here $\alpha$ denotes the liquid volume fraction that plays the role of the crowd density, $u$ denotes the velocity and $\pi\geq0$ denotes some singular pressure term appearing only when $\alpha=1$.\\
System (\ref{s0}-\ref{cs0}) is known in the literature as the {\it pressureless gas system with unilateral
constraint} and has been studied for example by {\sc F. Berthelin} in \cite{Be}. It can be interpreted as a coupling of two systems in the respective domains where $\alpha<1$ (liquid-gas mixture) and where $\alpha=1$ (pure liquid). 
{This system can be formally derived from bi-fluid system  (see f.i. {\sc F. Bouchut} {\it et al.} in \cite{BBCR} where a hierarchy for gas-liquid two-phase flows is also presented} ).
%
%
%
%
%

\bigskip

\noindent A generalization of  (\ref{s0}-\ref{cs0}) to the multi-dimensional viscous case is described {by} the compressible/incompressible Navier-Stokes type of system
\begin{equation}\label{CM0}
\left. \begin{array}{r}
 \partial_t\rho + \Div (\rho \vu) = 0\\
 0\leq\rho\leq \rho^*\\
 \partial_t (\rho \vu) + \Div (\rho \vu \otimes \vu) + \na p +\rho^*\na\pi-\Div\,\bf{S} = \vc{0}\\
  \Div\lr{\rho^*\vu}=0 \quad \text{a. e. in}\quad \{\rho=\rho^*\}\\  
  \pi \geq 0 \quad \text{a. e. in}\quad \{\rho=\rho^*\}\\
 \pi=0 \quad \text{a. e. in}\quad \{\rho<\rho^*\}\\
 \end{array}\right\} \quad \text{in} \ (0,T)\times\Omega,
\end{equation}
in which the homogeneous congestion constraint $\alpha{\leq} 1$ has been replaced by the inhomogeneous one $\rho\leq \rs(x)$.\\
The unknowns here are the density $\rho$, the velocity vector field $\vu$ and the pressure $\pi$, which is the Lagrange multiplier associated with the incompressibility constraint $\Div\lr{\rho^*\vu}=0$ a.e. in $\{\rho=\rho^*\}$. Note that as in the previous example, $\pi$ is apparent only in the congested regions $\{\rho=\rho^*\}$. In fact, conditions \eqref{CM0}$_{5}$, \eqref{CM0}$_{6}$ can be rewritten as one constraint
\eq{\rho\pi=\rho^*\pi\geq 0.\label{CON}}
Furthermore, the stress tensor $\vc{S}$ and the internal pressure $p$ are known functions of $\rho$ and $\vu$ which are typical for barotropic flow of Newtonian fluid, i.e.
\[ \vc{S}(\vu) = 2 \mu D(\vu) + \lambda \Div \vu \, \vc{I}, \quad \mu>0,\quad 2\mu+\lambda>0,\]
\eqh{p(\rho)=\rho^\gamma,\quad \gamma>1.}
System \eqref{CM0} mixes the free flow and congested regions denoted  by $\Omega_f(t)$ and $\Omega_c(t)$ respectively, where
$$\Omega_f(t) = \{x: \rho(t,x) <\rho^*(x)\}, \qquad 
     \Omega_c (t) = \{x: \rho(t,x)=\rho^*(x)\}.$$
It can be  thus seen as a free boundary problem, in which the the free interface separating $\Omega_c$ from $\Omega_f$ can become very irregular, as it is illustrated on the figure below. 

\bigskip
\begin{center}
\includegraphics{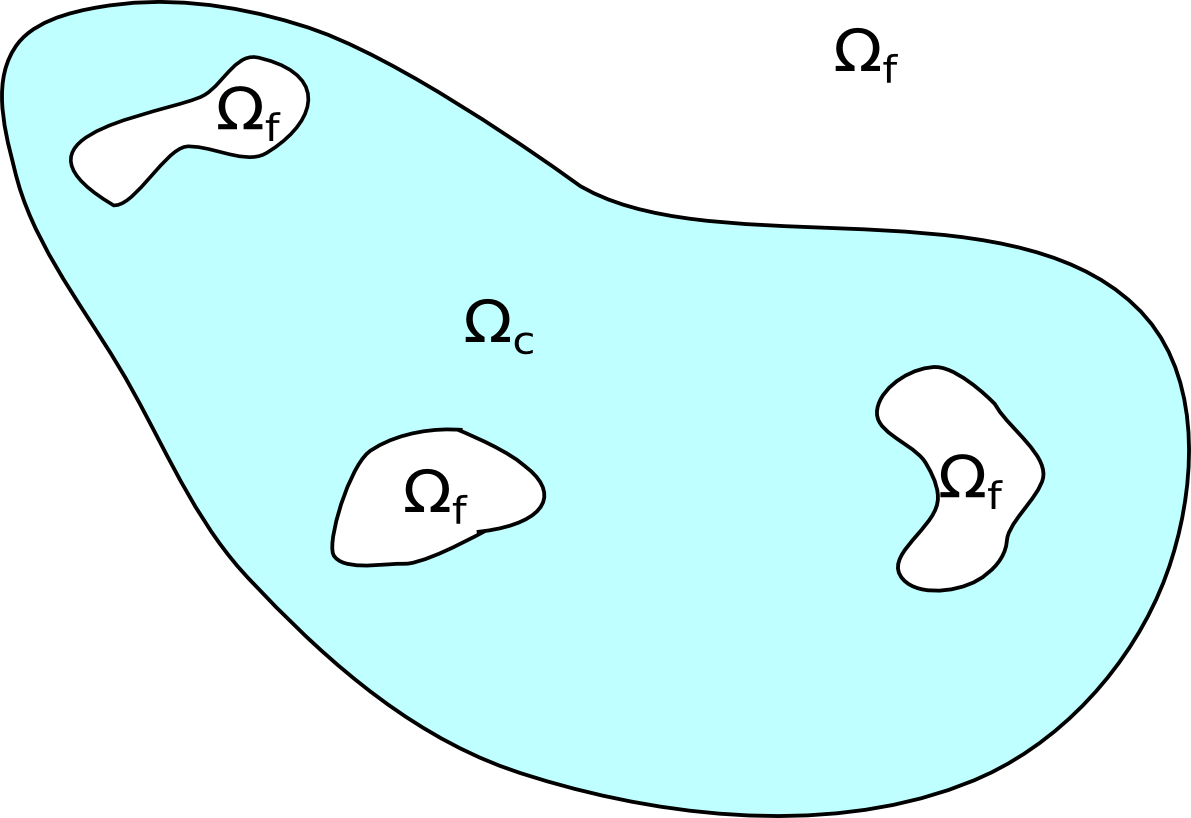}
\end{center}

\bigskip
The objective of this paper is to mathematically justify that the solution to problem \eqref{CM0} can be obtained as a  limit of $(\rho_n,\vu_n)$-- the solutions to the isentropic compressible Navier-Stokes equations 
\begin{equation}\label{CM0prim}
\left. \begin{array}{r}
 \partial_t\rho_n + \Div (\rho_n \vu_n) = 0\\
 \partial_t (\rho_n \vu_n) + \Div (\rho_n \vu_n \otimes \vu_n) + \na p +\rho^*\na\pi_{{n}}-\Div\,\bf{S} = \vc{0}
 \end{array}\right\} 
\end{equation}
where $\pi_{n}$ is some approximation of the limit pressure $\pi$. 
We want to find such an approximation that would guarantee uniform boundedness of the sequence approximating the density $\rho_{n}\leq\rho^*$. This  feature is very important for numerical purposes, see for example \cite{Ma}.

\bigskip

\noindent Below we present two possible ways of approximating the pressure $\pi$ appearing in the limit system \eqref{CM0}.

\medskip

\noindent{\bf The barotropic pressure.} The first kind of approximation uses the classical barotropic pressure
\begin{equation}\label{power_pres}\displaystyle \pi_{n}=\pi_{\gamma_n}=a\left(\frac{\rho_n}{\rho^*}\right)^{\gamma_n}
\end{equation}
with  $a >0$ fixed and the adiabatic exponent $\gamma_n>1$ being the approximation parameter.\\ 
Mathematical analysis of system \eqref{CM0prim} with the above pressure for fixed $\gamma_n$ is based on the existence theory for barotropic Navier-Stokes equations developed by {\sc P.--L. Lions} in \cite{PLL} and {\sc E. Feireisl} in \cite{EF2001}.\\
In this framework, only the homogeneous case $\rs=1$ has been studied. It has been justified by {\sc P.--L. Lions} and {\sc N.~Masmoudi} in \cite{LiMa} that the limit system \eqref{CM0} may be recovered from (\ref{CM0prim}-\ref{power_pres}) letting $\gamma_n\to\infty$. Later on, {\sc S.~Labb\'e} and {\sc E.~Maitre} performed the same limit passage for more complex system. They considered the viscosity coefficients $\mu,\ \lambda$ depending on the density and an additional surface tension term in the momentum equation \eqref{CM0}$_3$.\\
A similar asymptotic limit has been studied also  for a model of tumour growth by {\sc B.~Perthame},  {\sc F.~Guir\'os} and {\sc J.L.~V\'azquez} \cite{PeQuVa}. They used the cell population density model with the pressure of the form $\frac{m}{m-1}\lr{\rho/\rho^*}^{m-1}$ with parameter $m>1$ and the maximum packing density for the cells denoted by $\rho^*$ which is constant. In the limit $m\to\infty$ they obtained a free boundary model of Hele-Shaw type.

\bigskip

\noindent{\bf The singular pressure.} The second approximation uses a pressure that becomes singular close to some threshold value of the density $\rho^*$, for example
\begin{equation}\label{sing_pres}
 \displaystyle \pi_{\vep_n}= \vep_n\dfrac{\left(\dfrac{\rho_n}{\rho^*}\right)^\alpha}{\left(1-\dfrac{\rho_n}{\rho^*}\right)^\beta} 
 \end{equation} 
 with $\beta,\alpha>0$  fixed and  $\vep_n>0$ being the approximation parameter. \\
Such type of degeneration of the pressure can be used to model various phenomena. It appears in kinetic theory of dense gases where the interaction between the molecules is strongly repulsive at very short distance. The mutual reluctance of neighbouring molecules to share a certain amount of space (covolume) leads to the Van der Waals equation of state (see \cite{ChCo})
\[\pi = \frac{NRT}{V-Nb} - a\frac{N^2}{V^2}, \]
where $N$ is the number of molecules, $V$ the volume. The two terms on the right hand side (r.h.s.) represent respectively the repulsive and attractive forces. We see that the pressure becomes singular when $Nb$ approaches $V$, which corresponds to the state where motion is no longer possible. Other equations of state, modifying the representation of repulsive forces, were proposed for instance by {\sc N.F.~Carnahan} and {\sc K.E.~Starling} in \cite{CaSt}.  \\
Similar  form of the pressure has been also recently considered  by {\sc P.~Degond}, {\sc J.~Hua} and {\sc L.~Navoret} to model collective motion in \cite{DeHuNa}, and by {\sc F.~Berthelin} and {\sc D.~Broizat} to model traffic flow in \cite{BeBr}. Their approximation  is of the form $\nabla \pi_{\vep_n}$ with
 \[
 \displaystyle \pi_{\vep_n}= \dfrac{\vep_n}{\left(\dfrac{1}{\rho_n}-\dfrac{1}{\rho^*}\right)^\beta} ,
 \]
 while in \eqref{CM0prim} we take $\rs\nabla\pi_{\vep_n}$. However, as we will see later on, the factor $\rho^*$ is necessary in order to obtain the energy equality when $\rho^*$ is non-constant.\\
  Finally, systems involving such kind of degeneration are used in the theory of granular flows (see \cite{AnFoPo} and \cite{PiLe}). For instance, {\sc F.M.~Auzerais}, {\sc R.~Jackson} and {\sc W.B.~Russel} proposed in \cite{AuJaRu} a model for sedimentation using an empirical pressure of the form
\[\pi=\dfrac{C_0\phi^s}{\phi^*-\phi}\]
with $2\leq s\leq5$. In this framework the density of the fluid is replaced by $\phi$--the volume fraction of the solid phase ($0\leq\phi\leq1$) with some constant threshold value $\phi^*=0.64$.

\medskip

\noindent From the mathematical point of view, the first result  for system \eqref{CM0prim} with singular pressure is due to  {\sc E. Feireisl, H. Petzeltov\'a, E.~Rocca} and {\sc G.~Schimperna}, \cite{Fe}. They studied a model of two-phase compressible fluid flow with a Cahn-Hilliard type equation for a phase variable.  As a corollary of their result we get the existence of global in time weak solution for  \eqref{CM0prim}  and \eqref{sing_pres} with $\rho^*=const.$, $\beta>3$ and $\vep_n>0$ being fixed.

\medskip

\noindent{To our knowledge, the justification of the limit passage $\vep_n\to0$, which formally gives system \eqref{CM0}, is still an open problem in the general case.} The only result so far concerns  the one-dimensional homogeneous ($\rho^*=1$) case, studied by {\sc D. Bresch}, {\sc C. Perrin} and {\sc E. Zatorska} in \cite{BrPeZa}. First of all they proved that for $\beta,\gamma>1$ and $\vep>0$ fixed
there exists a {regular} solution $(\rho_\vep,u_\vep)$ to 
\eq{
 & \partial_t\rho_\vep +  \partial_x(\rho_\vep u_\vep) = 0, \\
 & \partial_t (\rho_\vep u_\vep) + \partial_x (\rho_\vep u_\vep \otimes u_\vep) + \partial_x p(\rho_\vep) + \vep\partial_x \frac{\rho_\vep^\gamma}{(1-\rho_\vep)^\beta} -(\lambda+\mu) \partial^2_{xx} u_\vep = 0,
\label{1D}}
such that 
\begin{equation}\label{b_rho}
0< c \le \rho_\varepsilon \le C(\varepsilon) < 1,
\end{equation}
for some constants $c$ and $C(\varepsilon)$.
Secondly, they justified that system \eqref{CM0} possesses a weak solution being a limit of regular solutions to \eqref{1D}.

\bigskip

\noindent The main difficulty in justifying the limit passage in both approaches is that the energy estimate does not give a uniform bound on the pressure term. More precisely, when $\rs=1$, the energy estimate for the power law pressure reads
 \[\Dt\intOB{\frac{1}{2}\rho_n |\vu_n|^2 + \frac{a}{\gamma_n-1}\rho_n^{\gamma_n}} + \intOB{\mu|\nabla \vu_n|^2+(\lambda+2\mu)(\Div \vu_n)^2} = 0  \]
 but it does not imply that $\rho^{\gamma_n}$ is bounded uniformly with respect to $\gamma_n\to\infty$.\\ 
 Similarly, for the degenerate pressure we have
 \[\Dt\intOB{\frac{1}{2}\rho_n |\vu_n|^2 + \rho_n\Gamma_{\vep_n}(\rho_n)} + \intOB{\mu|\nabla \vu_n|^2+(\lambda+2\mu)(\Div \vu_n)^2} = 0  \]
 with
 \[\Gamma_{\vep_n}(\rho_n)=\int_0^{\rho_n}{\dfrac{\pi_{\vep_n}(s)}{s^2}\mathrm{d}s},\]
but this does not yield the uniform pressure estimate either. Indeed, take for instance $\pi_{\vep_n}=\vep_n\rho_n^2(1-\rho_n)^{-4}$ for which $\Gamma_{\vep_n}=\vep_n/(3 (1-\rho_n)^{3})$, then the growth of $\rho_n\Gamma_{\vep_n}$ around singularity is one order less than of $\pi_{\vep_n}$. Additional information on the pressure is thus necessary and requires more sophisticated tools such as application of the Bogovskii operator.
%

\bigskip

\noindent The most important difference in these two ways of approximation lies in the uniform estimate of the density.  Indeed, approximation $\pi_{\gamma_n}=\rho_n^{\gamma_n}$ does not guarantee the validity of the congestion constraint $0\le \rho_{n}\le 1$ for fixed $\gamma_n$. The main advantage of approximation based on \eqref{sing_pres} is validity of this restriction uniformly with respect to $\vep_n$. In fact, the degenerate pressure plays the role of natural barrier (see \cite{Ma} by {\sc B. Maury}) and 
makes the second approach more suitable for numerical schemes. \\
\noindent  Our goal is therefore to extend the result from the previous work \cite{BrPeZa} to the global weak solutions framework in the multi-dimensional space case.
  The main difference between the one-dimensional case and the multi-dimensional case, is that in the latter, {\it  the sufficiently regular solutions are not known to exist}. Therefore, the strong convergence of the density and validity of the congestion constraint \eqref{CON} cannot be deduced directly from the {\it a priori} estimates. Eventually, the validity of the r.h.s. inequality in \eqref{b_rho} follows from an additional level of approximation using truncations of singular part of the pressure. 
  However, recovering system \eqref{CM0} after this step  requires {\it equi-integrability of the pressure term}, for which some restriction of the strength of singularity need to be imposed. The equi-integability of the pressure similar to \eqref{sing_pres} for $\rs=const.$ and $\beta>3 $ was proved by {\sc E.Feireisl} {\it et al.} in \cite{Fe} and we heavily relay on their approach to this issue. 

 The purpose of the paper is also to generalize \cite{BrPeZa} to the heterogeneous case, i.e. when the constant upper bound on the density is replaced by a prescribed function $\rho^*=\rho^*(x)>0$. This generalization has many applications. For instance, in  \cite{BeBr}, {\sc F.~Berthelin} and {\sc D.~Broizat} use a non-constant maximal constraint to study the dynamics of traffic jams and the influence of the number of lanes on the road.
It is also important in the study of models of flow through the closed pipes of non-uniform height $h^*(x)$ as said by {\sc F.~Berthelin} in \cite{Be}. In these models, the surface of the flowing fluid described by $h$ may  be either free when $h < h^*$ or additionally pressurized when it touches the pipe wall $h=h^*$, see the picture below.
\begin{center}
\includegraphics{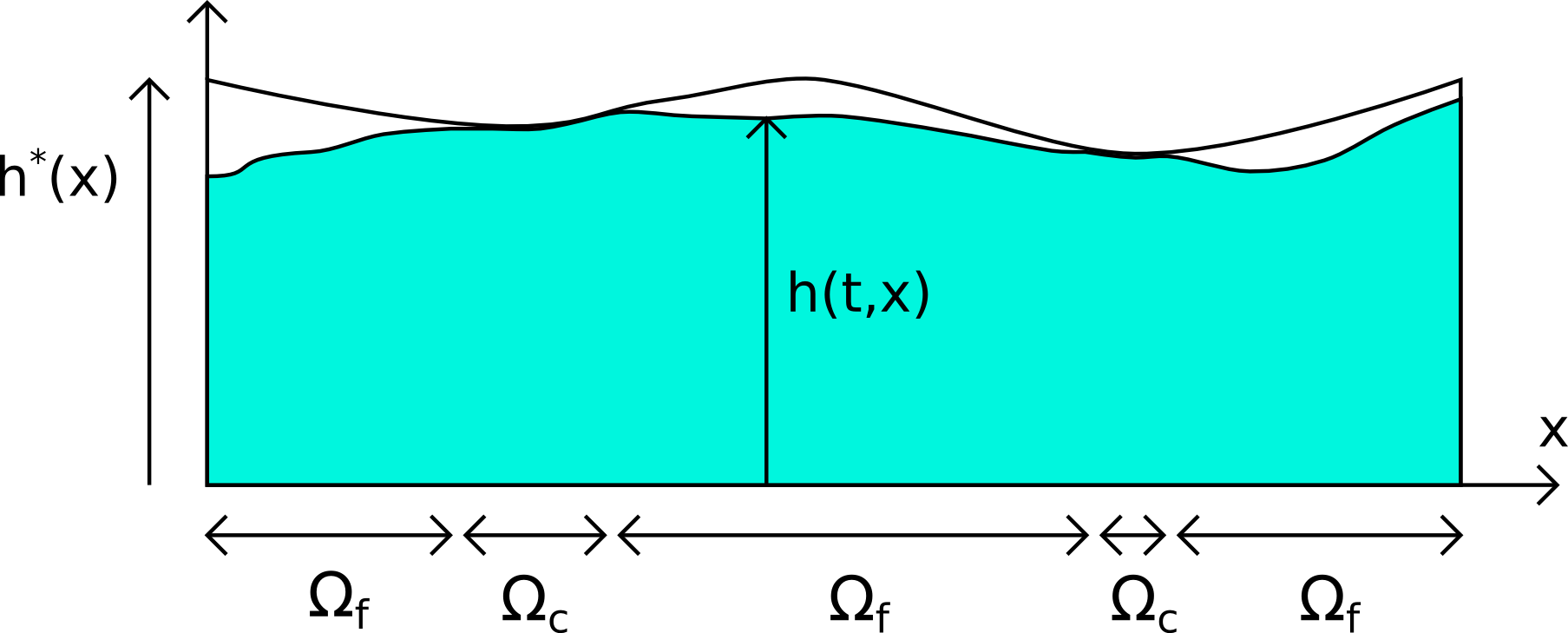}
\end{center}
The study of the non-viscous systems can be found in {\sc C.~Bourdarias, M.~Ersoy, S.~Gerbi}  \cite{BoEr}, see also references therein.

\section{Formulation of the main problem}
System \eqref{CM0} is supplemented with the homogeneous Dirichlet boundary conditions
\eq{\vu\big|_{\bom}=\vc{0}.
\label{bc}}
We assume that the threshold density $\rho^*=\rho^*(x)>0$ is a $\mathcal{C}^1(\overline{\Omega})$ function
 and that the initial data
\eq{\rho\big|_{t=0}=\rho_0,\quad (\rho\vu)\big|_{t=0}=\vc{m}_0\label{ic}}
satisfy
\begin{equation}\label{ini_data}
\begin{gathered}
 0\leq \rho_0(x) < \rho^*(x) \quad\text{a.e. in}\quad \Omega,\quad p(\rho_0) \in L^1(\Omega), \\
 \vc{m}_0\in (L^2(\Omega))^3, \quad \vc{m}_0\vc{1}_{\{\rho_0=0\}}=0 \quad\text{a.e. in}\quad \Omega,\\
 \frac{|\vc{m}_0|^2}{\rho_0}\vc{1}_{\{\rho_0>0\}} \in L^1(\Omega),
\end{gathered}
\end{equation}
\begin{equation}\label{add_hyp_ini}
M_0=\frac{1}{|\Omega|}\intO{\rho_0} <\inf_{x\in \Omega} \rs(x).
\end{equation}
\begin{rmk} Condition \eqref{add_hyp_ini} plays an essential role in the proof of uniform  $L^1((0,T)\times\Omega)$ bound for a sequence approximating $\pi$. Note that in the case $\rho^*$ constant, condition \eqref{add_hyp_ini} is directly satisfied since $\rho_0<\rho^*$.
\end{rmk}
Below we introduce the notion of a weak solution to system \eqref{CM0}.

\begin{df}[Weak solution of the limit system]\label{Def1}
A triple $(\rho,\vu,\pi)$ is called a weak solution to \eqref{CM0} with \eqref{bc} and \eqref{ic} if equations
\begin{equation*}
\begin{array}{c}
\partial_t\rho + \Div (\rho \vu) = 0,\\
 \partial_t (\rho \vu) + \Div (\rho \vu \otimes \vu) + \na p +\rho^*\na\pi-\Div\bf{S} = \vc{0}
\end{array}
\end{equation*}
are satisfied in the sense of distributions, the divergence free condition $\Div\lr{\rho^*\vu}=0$ is satisfied a.e. in $\{\rho=\rho^*\}$, the constraint $0\leq \rho\leq\rho^*$ is satisfied a.e. in $(0,T)\times \Omega$,
and the following regularity properties hold
\begin{align*}
&\rho\in \mathcal{C}([0,T];L^p(\Omega)), \quad  1\leq p < \infty, & \\
& \vu \in L^2(0,T;(W^{1,2}_0(\Omega))^3),\quad \dfrac{|\vc{m}|^2}{\rho} \in L^{\infty}(0,T; L^1(\Omega)),&\\
&\pi\in {\cal M}^+ ((0,T)\times \Omega).&
\end{align*}
Moreover, $\pi$ is sufficiently regular so that the condition
\eq{(\rho-\rho^*)\pi=0.\label{Cons}}
is satisfied in the sense of distributions.
\end{df}

\begin{rmk} Similarly to the homogeneous case studied bt {\sc P.-L. Lions} and {\sc N. Masmoudi} in {\rm \cite{LiMa}}, we can prove that the constraint $\eqref{CM0}_2$ and the divergence free condition $\eqref{CM0}_4$ are "compatible" (see {\rm Lemma \ref{LMP}}). More precisely, if $\vu$ belongs to $L^2(0,T;(W_0^{1,2}(\Omega))^3)$, $\rho$ belongs to $L^2((0,T)\times\Omega)$ and the couple $(\rho,\vu)$ satisfies the continuity equation $\eqref{CM0}_1$, then $\eqref{CM0}_2$ and $\eqref{CM0}_4$ with $0\leq \rho_0\leq \rho^*$ are equivalent.
 As it will be explained later, this fact is a natural consequence of the renormalized theory applied to the equations satisfied by the quantities $\rho/\rho^*$ and $\rho-\rho^*$.
\end{rmk}

\bigskip
\noindent Now we define  the notion of weak solution to approximate system \eqref{CM0prim} with the approximate pressure of the form
\begin{equation*}
\displaystyle \pi_\vep(r)= \left\{\begin{array}{ll} 
 \vep\dfrac{r^{\alpha}}{(1-r)^{\beta}} &\text{if } r < 1 \\
 \infty &\text{if } r \geq 1 \end{array}\right.\quad \alpha,\beta >3.
\end{equation*} 
 
\begin{df}[Weak solution of the approximate system]\label{weak_app_sol} A pair $(\rho_\vep,\vu_\vep)$ is called a weak solution to (\ref{CM0prim}) if it satisfies
\begin{itemize}
\item the approximate continuity equation
\eq{ \intTO{\rho_\vep\partial_t\phi} +\intO{\rho_0\phi(0)} + \intTO{\rho_\vep \vu_\vep \cdot \na \phi} \\
= \intO{(\rho_\vep \phi)(T) }\label{AC}}
for all $\phi\in \mathcal{D}([0,T]\times \Omega)$;
\item the approximate momentum equation
\eq{
 & \intTO{\rho_\vep \vu_\vep \cdot \partial_t \vp} + \intO{\vc{m}_0\cdot \vp(0)} \\
 +&\intTO{\rho_\vep(\vu_\vep\otimes\vu_\vep): \na \vp}- \intTO{\vc{S} : \na \vp}\\
 +&\intTO{p\Div\vp}
 +\intTO{\pi_\vep \left(\frac{\rho_\vep}{\rs}\right)\Div (\rs\vp)} \\
 &\qquad\qquad\qquad\qquad\qquad\qquad\qquad= \intO{ (\vc{m}_\vep \cdot\vp)(T)},
\label{AM}}
for all $\vp\in (\mathcal{D}([0,T]\times \Omega))^3$.
 \end{itemize}
 \end{df}
 \bigskip
 
\noindent The objective of this paper is to  justify that the weak solution from \mbox{Definition \ref{Def1}} can be obtained as a limit of weak solutions from Definition \ref{weak_app_sol}. \\

\medskip

\noindent The main theorem of this paper reads.
\begin{Theorem} \label{main}
\begin{enumerate}
 \item Let $\vep>0$ be fixed, then there exists a global weak solution $(\rho_\vep,\vu_\vep)$ to \eqref{CM0prim} in the sense of Definition \ref{weak_app_sol}, moreover,
\eq{0\leq\re\leq\rs.\label{un_main}}
 \item For $\vep\to0$, there exists a subsequence $(\re,\ue,\pi_{\vep})$ converging to  $(\rho,\vu,\pi)$ a solution of system \eqref{CM0} in the sense of Definition \ref{Def1}.
More precisely
\eqh{
&\re \to \rho \quad\text{weakly\ in }\quad L^p((0,T)\times\Omega) \quad \forall 1\leq p <+\infty,\\
&\ue \to \vu \quad\text{weakly\ in }\quad L^2(0,T;(W^{1,2}_0(\Omega))^3),\\
&\pi_\vep \to \pi \quad\text{weakly\ in }\quad {\cal M}^+((0,T)\times\Omega).
}
\end{enumerate}
\end{Theorem}

\medskip

The paper is organized as follows. In Section \ref{S_rho1}
we present details of approximation and prove the first part of Theorem \ref{main}. Then, in Section \ref{S_ep}, we recover the original system by letting $\vep\to 0$.  Section \ref{Appendix} is the Appendix in which we recall some basic facts about the Bogovskii operator, the Riesz transform and the renormalized continuity equation that are used in several places in the course of the proof.

 \section{Existence of approximate solutions}\label{S_rho1}

In order to prove the first part of Theorem \ref{main}, we consider further approximation with additional truncation parameter $\delta$ and artificial pressure $\kappa\rho^{K}$ where $K$ is sufficiently large positive number that will be determined later on. Then we will show that  (\ref{AC}-\ref{AM}) can be recovered by letting $\delta\to0$ and $\kappa\to0$ respectively. In  this section $\vep$ is fixed and we drop it when no confusion can arise.

\subsection{Basic level of approximation}

For fixed $\kappa,\delta>0$, we consider the basic level of approximation 

%


\eq{\label{appr_model_delta}
 & \partial_t\rho_\delta + \Div (\rho_\delta \vu_\delta) = 0 \\
 & \partial_t (\rho_\delta \vu_\delta) + \Div (\rho_\delta \vu_\delta \otimes \vu_\delta) +  \na p(\rho_\delta) +\rs\na\pi_{\kappa,\delta}\lr{\frac{\rho_\delta}{\rs}} -\Div\vc{S} = 0
}
with the approximate pressure 
given by
\begin{equation}\label{pressure_ap_p2}
\displaystyle \pi_{\kappa,\delta}(s) =\kappa s^K+\left\{
\begin{array}{ll}
\vspace{0.2cm}
\displaystyle\vep \frac{s^{\alpha}}{(1-s)^{\beta}} & \text{ if } s < 1-\delta, \\
\displaystyle\vep \frac{s^{\alpha}}{\delta^{\beta}} & \text{ if } s\geq 1-\delta .            
\end{array}\right.
\end{equation}
%
%
%
%

\bigskip

For all parameters fixed, the approximate pressure $\pi_{\kappa,\delta}$ is a monotone increasing function of $\rho$. For such pressure the issue of existence of global in time weak solutions is a straightforward adaptation of the proof from the case of barotropic system see for example \cite{PLL}, \cite{FN} or \cite{NS}. Below we will only make some general comments concerning the {\it a priori} estimates that are necessary to state the analogous existence result.

First of all, let us formally write the basic energy equality. Multiplying the second equation of \eqref{appr_model_delta} by $\vu$, integrating by parts and using the mass equation we easily get
\[\Dt \intO{\frac{1}{2} \rho_\delta |\ud|^2} + \intO{\big(\nabla p+ \rho^*\nabla \pi_{\kappa,\delta}\big)\cdot \ud}+ \intO{\vc{S}:\nabla \vu_\delta}=0.\]
Let us explain how to deal with the term coming from the pressure for reader's convenience but also to check that everything works with the heterogeneous maximal density:
\begin{align*}
\intO{\big(\nabla p+ \rho^*\nabla \pi_{\kappa,\delta}\big)\cdot \ud} & =  \intO{\gamma\rho_\delta^{\gamma-2}\nabla\rho_\delta\cdot (\rho_\delta\ud)}\\
& \quad + \intO{\dfrac{\rho^*}{\rho_\delta}\pi_{\kappa,\delta}'\left(\dfrac{\rho_\delta}{\rho^*}\right)\nabla\left(\dfrac{\rho_\delta}{\rho^*}\right) \cdot(\rho_\delta\ud)} \\
& = \intO{\nabla\lr{\dfrac{\gamma}{\gamma-1}\rho_\delta^{\gamma-1}+Q_{\kappa,\delta}\lr{\dfrac{\rho_\delta}{\rho^*}}} \cdot \big(\rho_\delta\vu_\delta\big)} \\
& = -\intO{\lr{\dfrac{\gamma}{\gamma-1}\rho_\delta^{\gamma-1}+Q_{\kappa,\delta}\lr{\dfrac{\rho_\delta}{\rho^*}}} \Div\,(\rho_\delta\vu_\delta)} \\
& = \intO{\lr{\dfrac{\gamma}{\gamma-1}\rho_\delta^{\gamma-1}+Q_{\kappa,\delta}\lr{\dfrac{\rho_\delta}{\rho^*}}}  \partial_t\rho_\delta}
\end{align*}
where we denoted $\displaystyle Q_{\kappa,\delta}'(r)=\frac{\pi_{\kappa,\delta}'(r)}{r}$. Therefore, the energy equality reads 
\eq{ \label{energ_est}
\Dt &\intO{\lr{\frac{1}{2} \rho_\delta |\vu_\delta|^2 + \frac{1}{\gamma-1}\rho_\delta^\gamma+ \rd\Gamma_{\kappa,\delta}\lr{\frac{\rho_\delta}{\rho^*}}}} \\
&+ \intO{\vc{S}:\nabla\vu_\delta}=0
}
with $\Gamma_{\kappa,\delta}$ such that $\Gamma_{\kappa,\delta}(r)+r\Gamma_{\kappa,\delta}'(r)=Q_{\kappa,\delta}(r)$. To find the expression of $\Gamma_{\kappa,\delta}$ in terms of $\pi_{\kappa,\delta}$, we use the definition of $Q_{\kappa,\delta}$, integrate by part and use that
$[\pi_{\kappa,\delta}(r)/r]\vert_{r=0} =0$. We get
\eq{\label{Gamma}
\Gamma_{\kappa,\delta}(s) =  \int_0^{s} {\frac{\pi_{\kappa,\delta}(r)}{r^2}\, {\rm d}r}.}
Integrating \eqref{energ_est} with respect to time gives rise to the following estimates
\begin{equation} \label{un_delta}
\begin{gathered}
\sup_{t\in[0,T]}\lr{\|\sqrt{\rho_\delta} \vu_\delta (t)\|_{L^2(\Omega)} 
+\|\rho_\delta (t)\|_{L^{\gamma}(\Omega)} 
+\left\|\rho_\delta \Gamma_{\kappa,\delta}\lr{\frac{\rho_\delta}{\rho^*}}(t) \right\|_{L^1(\Omega)} }\leq C,\\
\intT{{\|\vu_\delta\|_{W^{1,2}(\Omega)}^2}} \leq C.
\end{gathered}
\end{equation}
In particular, since 
\begin{align*}
\displaystyle \rho_\delta \Gamma_{\kappa,\delta}\lr{\frac{\rho_\delta}{\rho^*}}= \rho_\delta\int_0^{\rho_\delta}{\dfrac{\pi_{\kappa,\delta}(s)}{s^2}\,\mathrm{d}s}
\geq 
\dfrac{\kappa}{K-1}\rho_\delta^K
\end{align*} 
we obtain that $\rho_\delta$ is bounded in $L^\infty(0,T; L^K(\Omega))$.
\begin{rmk} This computation shows that the approximate of the form
\[\displaystyle\vep\rho^* \nabla \dfrac{\left(\dfrac{\rho}{\rho^*}\right)^\alpha}{\left(1-\dfrac{\rho}{\rho^*}\right)^\beta}\]
gives a good contribution to the energy. The form of the so-called potential energy from \eqref{energ_est} is a natural extension of the one obtained in {\rm \cite{BrPeZa}} and {\rm \cite{Fe}} to the heterogeneus case.
 This  would not be the case for the approximation 
\[\displaystyle\vep\nabla \dfrac{1}{\left(\dfrac{1}{\rho}-\dfrac{1}{\rho^*}\right)^\beta}\] 
proposed in {\rm \cite{BeBr}} or {\rm \cite{DeHuNa}}.
\end{rmk}
 The above estimate may be used to improve integrability of $\rho_\delta$. Indeed, taking $K$ sufficiently large, say $K>4$, we can test \eqref{AM} by
  \[\displaystyle \varphi=\dfrac{\psi(t)}{\rho^*}{\cal B}\lr{\rho_{\delta}-\frac{1}{|\Omega|}\int_\Omega{\rho_{\delta}(y){\rm d}y}},\]
  where ${\cal B}$ is the Bogovskii operator (for definition and properties see Lemma \ref{lem_bog} and Proposition \ref{prop_bog} in Appendix).    
The details of this testing will be given in Section \ref{uniform_est}. In particular, it gives an additional control (but not uniform with respect to $\kappa$) for the density
\eq{\|\rho_{\delta}\|_{L^{K+1}((0,T)\times \Omega)}\leq C.\label{un_bog}} 
Combination of these estimates can be used to deduce that
\eq{\|\rho_\delta|\vu_\delta|^2\|_{L^2(0,T;L^{\frac{6K}{4K+3}}(\Omega))}
+\|\rho_\delta\vu_\delta\|_{L^\infty(0,T; L^{\frac{2K}{K + 1}}(\Omega))\cap L^2(0,T; L^{\frac{6K}{K + 6}}(\Omega))}\leq C.\label{interp}}
These observations allow us to apply methods developed for the barotropic system from \cite{FN, PLL} to justify the following statement.
\begin{prop}[Existence of weak solutions] Let $\vep$,  $\kappa$, $\delta$ be fixed and positive. Then, there exists a couple $(\rho_{\delta}, \vu_{\delta})$ solving (\ref{appr_model_delta}-\ref{pressure_ap_p2}) in the sense of distributions on $(0,T)\times\Omega$ with the following regularity properties
\begin{equation}\label{0_reg}
\begin{gathered}
\rho_{\delta} \in L^\infty(0,T;L^{K}(\Omega)) \cap L^{K+1}((0,T)\times \Omega),\\
\rho_{\delta} \in \mathcal{C}([0,T];L^{K}_{\rm weak}(\Omega)) \cap \mathcal{C}([0,T];L^p(\Omega)) , \quad 1\leq p <K, \\
\rho_{\delta} \geq 0 \quad\text{a.e in }  (0,T)\times\Omega, 
 \\
 \vu_{\delta}\in L^2(0,T;(W_0^{1,2}(\Omega))^3), 
\\
 \rho_{\delta} |\vu_{\delta}|^2 \in L^\infty(0,T;L^1(\Omega))\cap L^2(0,T;L^{\frac{6K}{4K+3}}(\Omega)),\\
 \rho_{\delta} \vu_{\delta} \in \mathcal{C}([0,T]; (L^{\frac{2K}{K + 1}}_{\rm weak}(\Omega))^ 3)\cap L^2(0,T;(L^{\frac{6K}{K+6}}(\Omega))^3)
\end{gathered}
\end{equation}
and such that the energy inequality
\eq{ \label{energ_in}
\Dt \intOB{\frac{1}{2} \rho_\delta |\vu_\delta|^2 + \frac{1}{\gamma-1}\rho_\delta^\gamma
+ \rd\Gamma_{\kappa,\delta}\lr{\frac{\rho_\delta}{\rho^*}}} 
+ \intOB{\vc{S}:\nabla\vu_\delta} \leq0 
}
 is satisfied in the sense of distributions with respect to time.
 
Moreover $(\rho_{\delta}, \vu_{\delta})$ extended by $0$ outside $\Omega$ is the renormalized solution to the continuity equation in the sense of Definition \ref{df2}.
\end{prop}

\subsection{Uniform estimates}\label{uniform_est}
The goal of this subsection is to provide estimates which are uniform with respect to $\delta$. One of them is estimate \eqref{un_delta}. Note, however, that it does not assure boundedness of the singular part of the pressure $\pi_{\kappa,\delta}$. In fact, such a bound follows from Bogovskii estimate announced in the previous section and giving rise to \eqref{un_bog}. We now check that this estimate assures also the uniform $L^1((0,T)\times\Omega)$ for the pressure.

 \bigskip

\ni{\bf $L^1$ bound of the pressure.} 
From estimates \eqref{un_delta} it follows that $\rd\in L^{K}((0,T)\times\Omega)$ with $K \geq 4$, therefore we can test the momentum equation by
\[\vp(t,x) = \frac{\psi(t)}{\rho^*}\mathcal{B}\left(\rd - \Ov{\rd}\right),\qquad  \psi(t)\in \mathcal{C}^\infty_0((0,T)),\quad \psi\geq0,\]
where we denoted $\Ov{\rd}= \frac{1}{|\Omega|}\int_\Omega{\rd}{\rm d}y$. 
We then obtain 

\eq{\label{B}
&\intTO{ \psi{\pi_{\kappa,\delta}\lr{\frac{\rd}{\rs}}\left(\rd - \overline{\rd}\right)}} 
+\intTO{\psi{ p(\rd)\frac{\rd}{\rs}}}\\
& \qquad\qquad=  \intTO{\partial_t(\rd\ud)\cdot \vp} -\intTO{\rd(\ud\otimes \ud) : \na \vp} \\
& \qquad\qquad\quad+ \intTO{\vc{S} : \na \vp}
+\intTO{\psi{ p(\rd)\frac{\Ov{\rd}}{\rs}}}\\
& \qquad\qquad\quad-\intTO{\psi{ p(\rd)\na\lr{\frac{1}{\rs}}\cdot{\cal B}(\rd-\Ov{\rd})}}\\
& \qquad\qquad=\sum_{i=1}^5 I_{i}}
A priori estimates obtained in \eqref{un_delta} allow to control the r.h.s. Indeed for the first term we may write
\begin{align*}
I_{1}&=-\intTO{\psi'{\frac{\rd}{\rho^*}\ud \cdot \mathcal{B} \left(\rd - \overline{\rd}\right)}} -\intTO{\psi {\frac{\rd}{\rho^*}\ud \cdot \partial_t\mathcal{B} \left(\rd - \overline{\rd}\right)}}.
\end{align*}
Therefore, if only $K>3$ and since $\psi, \psi'$ and $\lr{\rho^*}^{-1}$ are bounded in $L^\infty((0,T)\times\Omega)$ and thanks to the mass equation we get
\eqh{
|I_{1}|& \leq  C\intT{\|{\rd}\|_{L^{3}(\Omega)}\|\ud\|_{L^6(\Omega)}\|\mathcal{B} \left(\rd - \overline{\rd}\right)\|_{L^{6/5}(\Omega)}}\\
&  \quad+ C\intT{\|{\rd}\|_{L^3(\Omega)}\|\ud\|_{L^6(\Omega)}\|\mathcal{B} \left(\Div (\rd\ud )\right)\|_{L^{6/5}(\Omega)}} \\
& \leq  C\intT{\|{\rd}\|^2_{L^{3}(\Omega)}\|\ud\|_{L^6(\Omega)}}\\
&  \quad+ C\intT{\|{\rd}\|_{L^3(\Omega)}^2\|\ud\|_{L^6(\Omega)}^2} \\
& \leq C.
}
Similarly, for the second term, we have
\eqh{
I_2=&
\intTO{\psi\frac{\rd}{\rho^*}\lr{\vu_\delta\otimes\vu_\delta}:\nabla \mathcal{B}\left(\rd - \Ov{\rd}\right)}\\
&+\intTO{\psi{\rd(\ud\otimes \ud)\na\lr{\frac{1}{\rs}}\cdot{\cal B}(\re-\Ov{\re})}},}
thus
\begin{align*}
|I_{2}|&\leq C\intTO{{\rd|\ud|^2\lr{|\na \mathcal{B}\big(\rd-\overline{\rd}\big) |+| \mathcal{B}\big(\rd-\overline{\rd}\big) |}}}\\
&\leq C\intT{\|\rd\|^2_{L^3(\Omega)}\|\ud\|^2_{L^6(\Omega)}}  \leq C.
\end{align*}
For the stress tensor $I_3$ we may write
\eqh{
|I_{3}|& \leq C\intT{\|\na\ud\|_{L^2(\Omega)}\|\rd\|_{L^2(\Omega)}}\leq  C.
}
The remaining pressure terms $I_{4}$ and $I_5$ can be easily controlled using the uniform $L^\infty(0,T;L^K(\Omega))$ bound for $\rd$ obtained in \eqref{un_delta}, if only $K$ is sufficiently large
\eqh{
|I_4|\leq C\intT{\|\rd\|_{L^1(\Omega)}\|\rd\|_{L^\gamma(\Omega)}^\gamma}\leq C
}
for $K>\gamma$ and
\eqh{
|I_5|\leq C\intT{\|\rd\|_{L^{\frac{2\gamma K}{K-1}}(\Omega)}^\gamma\|\rd\ud\|_{L^{\frac{2K}{K+1}}(\Omega)}}
\leq C}
for $K>2\gamma+1$. Note that in above estimate we essentially use the $L^\infty((0,T)\times\Omega)$ bound on $\nabla\lr{\dfrac{1}{\rho^*}}$, which follows from assumptions on $\rho^*$ in Theorem \ref{main}.
Collecting these estimates, we verify that the l.h.s. of \eqref{B} is bounded uniformly with respect to $\delta$.
Now, let us split the first term into two parts as follows
\begin{align} \label{splitting}
&\intTO{\psi{\Pd\left(\rd -\Ov{\rd}\right)}} \nonumber \\
\displaystyle &\quad = \intTO{\psi{\Pd\left(\rd -\Ov{\rd}\right)\vc{1}_{\displaystyle\left\lbrace\rd < \dfrac{\inf\rs + M_0}{2}\right\rbrace}}} \\
\displaystyle & \qquad + \intTO{\psi{\Pd\left(\rd -\Ov{\rd}\right)\vc{1}_{\displaystyle\left\lbrace\rd \geq \dfrac{\inf\rs + M_0}{2}\right\rbrace}}}\nonumber
\end{align}
In the first integral, the term $\Pd$ is far from the singularity $\rs$, so the integral is finite. 
 For the second integral we may write
\begin{align*} \displaystyle&\intTO{\psi{\Pd\left(\rd -\Ov{\rd}\right)\vc{1}_{\displaystyle\left\lbrace\rd \geq \dfrac{\inf\rs + M_0}{2}\right\rbrace}}}\\
\displaystyle& \quad \geq \frac{\inf \rs -M_0}{2} \intTO{\psi{\Pd}\vc{1}_{\displaystyle\left\lbrace\rd \geq \dfrac{\inf\rs + M_0}{2}\right\rbrace}}.
\end{align*}
Thanks to \eqref{add_hyp_ini}, we deduce from the second integral of \eqref{splitting} that 
\eqh{\intTO{\Pd} \leq C
\label{pressure_bound_delta}}
and then the control of the first integral of \eqref{splitting} yields
\eqh{{\intTO{\rd \Pd}} \leq C.}
 In particular, as mentioned in \eqref{un_bog}, $\rd$ is bounded in  $L^{K+1}((0,T)\times\Omega)$
uniformly with respect to $\delta$, but not uniformly with respect to $\kappa$.

\bigskip

\ni{\bf Equi-integrability of the pressure.} In order to perform the limit passage $\delta\rightarrow 0$ and to recover system (\ref{AC}-\ref{AM}) with parameter $\kappa$, the weak limit of $\Pd$ has to be more regular than merely a positive measure. To show that it converges weakly in $L^1((0,T)\times \Omega)$ to $\pi_{\kappa}\lr{\dfrac{\rho_\kappa}{\rho^*}}$ we want to apply the {\sc De~La~Vall\'ee-Poussin} criterion and we test the momentum equation by
\[\varphi(t,x) = \frac{\psi(t)}{\rho^*}\mathcal{B}\left(\eta_\delta\left(\dfrac{\rd}{\rho^*}\right) - \Ov{\eta_\delta\left(\dfrac{\rd}{\rho^*}\right)}\right), \quad\psi(t)\in {\mathcal C}^\infty_0((0,T)),\quad \psi\geq0,\]
where
\begin{equation}\label{eta_function}
\eta_\delta(s) = \begin{cases}  \,\,\log (1-s) & \text{ if } s \leq 1-\delta, \\
                                 \,\,\log (\delta)  & \text{ if } s > 1-\delta. \end{cases}
\end{equation} 
This testing results in
\begin{align}\label{e-i}
&\intTO{ \psi\lr{\pi_{\delta}\lr{\frac{\rd}{\rs}}+\frac{p(\rd)}{\rs}}\left(\eta_\delta\left(\dfrac{\rd}{\rho^*}\right) - \overline{\eta_\delta\left(\dfrac{\rd}{\rho^*}\right)}\right)} \nonumber \\
& =  \intTO{\partial_t(\rd\ud)\cdot \vp} -\intTO{\re(\ue\otimes \ue) : \na \vp} \nonumber\\
& \quad+ \intTO{\vc{S} : \na \vp}-\intTO{\psi{ p(\rd)\na\lr{\frac{1}{\rs}}\cdot{\cal B}\lr{\eta_\delta\left(\dfrac{\rd}{\rho^*}\right)-\Ov{\eta_\delta\left(\dfrac{\rd}{\rho^*}\right)}}}}\nonumber\\
& =\sum_{i=1}^4 J_{i}.\end{align}
Similarly as in the previous paragraph, the most demanding term $J_1$ equals to
\eqh{J_1 =  &-\intTO{\psi' {\frac{\rd}{\rho^*}\ud \cdot \mathcal{B} \left(\eta_\delta\left(\dfrac{\rd}{\rho^*}\right) - \overline{\eta_\delta\left(\dfrac{\rd}{\rho^*}\right)}\right)}} \\
&-\intTO{\psi {\frac{\rd}{\rho^*}\ud \cdot \partial_t\mathcal{B} \left(\eta_\delta\left(\dfrac{\rd}{\rho^*}\right) - \overline{\eta_\delta\left(\dfrac{\rd}{\rho^*}\right)}\right)}}}
We can generalize the notion of renormalized solutions of the continuity equation (see the Appendix and \cite{NS} section 6.2) to functions 
\[b_k(s) = \begin{cases} & b(s) \text{ if } s\in[0,k) \\
                         & b(k) \text{ if } s\in[k,\infty) \end{cases} \quad \text{with } b\in {\cal C}^1[0,\infty)\]
 to deduce that  $\eta_\delta$ satisfies the equation
\begin{align} \label{renorm_eta}
&\partial_t \eta_\delta\left(\dfrac{\rd}{\rho^*}\right) + \Div\left(\eta_\delta\left(\dfrac{\rd}{\rho^*}\right) \ud\right) + \left((\eta_\delta)'_+\left(\dfrac{\rd}{\rho^*}\right)\dfrac{\rd}{\rho^*} -\eta_\delta\left(\dfrac{\rd}{\rho^*}\right)\right) \Div \ud \nonumber\\
&\quad + (\eta_\delta)'_+\left(\dfrac{\rd}{\rho^*}\right)\dfrac{\rd}{\rho^*}\ud\cdot\na\log\rs=0.
\end{align}
In the above formula $(\eta_\delta)'_+$ denotes the right derivative of $\eta_\delta$
\[(\eta_\delta)'_+ (s)=\begin{cases} \, \dfrac{1}{1-s} &\text{ if } (x,t)\in\{s(x,t) < 1-\delta\} \\
\, 0 & \text { if }(x,t)\in\{s(x,t) \geq 1-\delta\}\end{cases}.\]
We thus obtain
\eq{\label{J1}
J_1=& -\intTO{\psi'{\frac{\rd}{\rho^*}\ud \cdot \mathcal{B} \left(\eta_\delta\left(\dfrac{\rd}{\rho^*}\right) - \overline{\eta_\delta\left(\dfrac{\rd}{\rho^*}\right)}\right)}} \\
&- \intTO{\psi {\frac{\rd}{\rho^*}\ud \cdot \mathcal{B} \left(\Div\left(\eta_\delta\left(\dfrac{\rd}{\rho^*}\right)\ud\right)\right)}} \\
& -\intTO{ \psi\frac{\rd}{\rho^*}\ud\cdot\mathcal{B}\left[F_\delta\right]}
}
where we denoted
\begin{align*}
F_\delta=&\left(\eta_\delta'\left(\dfrac{\rd}{\rho^*}\right)\dfrac{\rd}{\rho^*} -\dfrac{\rd}{\rho^*}\right) \Div \ud + \eta_\delta'\left(\dfrac{\rd}{\rho^*}\right)\dfrac{\rd}{\rho^*}\ud\cdot\na\log\rs \\ &-\frac{1}{|\Omega|}\int_\Omega{\left[\left(\eta_\delta'\left(\dfrac{\rd}{\rho^*}\right)\dfrac{\rd}{\rho^*} -\eta_\delta\left(\dfrac{\rd}{\rho^*}\right)\right) \Div \ud + \eta_\delta'\left(\dfrac{\rd}{\rho^*}\right)\dfrac{\rd}{\rho^*}\ud\cdot\na\log\rs \right]{\rm d}y}.
\end{align*}
Therefore, to control the last integral of $J_1$ we need bounds on $\eta_\delta\left(\dfrac{\rd}{\rho^*}\right)$ and $\eta_\delta'\left(\dfrac{\rd}{\rho^*}\right)\dfrac{\rd}{\rho^*}$. To simplify, we assume that $\alpha, \beta$ are integers. For $\dfrac{\rd}{\rho^*} < 1-\delta$ we have
\begin{align*}
\displaystyle\rd \Gamma_{\kappa,\delta}\lr{\dfrac{\rd}{\rho^*}} & 
\geq \rd \int_0^{\frac{\rd}{\rho^*}}{\vep\dfrac{s^{\alpha-2}}{(1-s)^\beta} {\rm d}s} \\
                    \displaystyle    & = \vep \rd \int_{1-\frac{\rd}{\rho^*}}^1{\dfrac{(1-s)^{\alpha-2}}{s^\beta}ds} \\
                        & = \vep \rd \sum_{k=0}^{\alpha-2}\begin{pmatrix} \alpha-2 \\ k \end{pmatrix} (-1)^k\int_{1-\frac{\rd}{\rho^*}}^1{s^{k-\beta}ds}\\
                        & = \dfrac{\vep \rd}{\beta-1}\dfrac{1}{\lr{1-\dfrac{\rd}{\rho^*}}^{\beta-1}}-\dfrac{\vep \rd}{\beta-1}\\
                        &\quad+\vep \rd \sum_{k=1}^{\alpha-2}\begin{pmatrix} \alpha-2 \\ k \end{pmatrix} (-1)^k\int_{1-\frac{\rd}{\rho^*}}^1{s^{k-\beta}ds} \\
                        & \geq \dfrac{\vep \rd}{\lr{1-\dfrac{\rd}{\rho^*}}^{\beta-1}}\left(\dfrac{1}{\beta-1}-\delta C(\alpha,\beta,\delta)\right) -C_2
\end{align*}
where, in the last inequality, $C_2$ controls all the non-degenerate terms. The constant $C(\alpha,\beta,\delta)$ is such that $\lim_{\delta\rightarrow 0} |C| < +\infty$ and it controls all the degenerate  terms of lower order. Taking $\delta$ sufficiently small, we ensure that 
\[\delta|C(\alpha,\beta,\delta)|<\dfrac{1}{2(\beta-1)}.\]
Thus it follows that
\begin{equation}\label{ineg_equi}
\displaystyle\rd \Gamma_{\kappa,\delta}\lr{\frac{\rd}{\rho^*}}  \geq C_1\frac{\vep\rd}{\lr{1-\frac{\rd}{\rho^*}}^{\beta-1}} - C_2, \qquad C_1 >0 .
\end{equation}
Now, since for $\dfrac{\rd}{\rho^*} < 1-\delta$
\[\eta'_\delta\left(\dfrac{\rd}{\rho^*}\right) = \dfrac{1}{\left(1-\dfrac{\rd}{\rho^*}\right)}\]
we obtain using the Young inequality ($\beta>3$)
\[C_1|\rd\eta'_\delta(\rd)|^2 \leq \rd \Gamma_{\kappa,\delta}\lr{\frac{\rd}{\rho^*}} +C_2.\]
This inequality is still satisfied for $\dfrac{\rd}{\rho^*} \geq 1-\delta$ because $\eta'\equiv 0$ on this set.\\
Concerning the control of $\eta_\delta$, if $\dfrac{\rd}{\rho^*} < 1-\delta$ 
we observe that
\[(1-R)^{\beta-1}|\log(1-R)|^q \underset{R\rightarrow 1^-}{\longrightarrow} 0 \quad \forall 1\leq q< \infty\]
which proves that there exists $C_1$ and $C_2$ uniform with respect to $\delta$ such that
\begin{equation} \label{bound_eta}
C_1(q,\vep,\beta)\left|\eta_\delta\left(\dfrac{\rd}{\rho^*}\right)\right|^q \leq \rd \Gamma_{\kappa,\delta}\left(\dfrac{\rd}{\rho^*}\right) +C_2(q,\vep,\beta)\quad \forall q\in [1,\infty). \end{equation}
Since $\eta_\delta\left(\dfrac{\rd}{\rho^*}\right)$ is constant for $\dfrac{\rd}{\rho^*} \geq 1-\delta$, we deduce that this inequality is satisfied for all $\dfrac{\rd}{\rho^*}$.

These estimates allow to control $\mathcal{B}\left[F_\delta\right] $ from \eqref{J1}
uniformly in $L^2(0,T;L^q(\Omega))$ with $q<3/2$. So, the full integrant in the third term of $J_1$ is controlled provided $\rd\ud\in L^2(0,T;(L^{p}(\Omega))^3)$ with $p> 3$, which due to \eqref{interp} asks for $K>6$. 

%
The other $J_k$ from \eqref{e-i} are bounded thanks to \eqref{bound_eta} similarly as in the case of $L^1$ bound on the pressure obtained in previous paragraph.

Finally, from \eqref{e-i} we deduce that
\begin{equation}\label{equiintegr}
\intTO{\pi_{\delta}\lr{\frac{\rd}{\rs}}\eta_\delta\left(\dfrac{\rd}{\rho^*}\right)} \leq C, 
\end{equation} 
which implies, thanks to the {\sc De La Vall\'ee-Poussin} criterion, the equi-integrability of the pressure $\pi_{\delta}$.

\subsection{Passage to the limit $\delta \to 0$, $\kappa\to0$}
First, let us perform the limit passage $\delta\to0$. In fact the second limit passage $\kappa\to0$ is an immediate consequence of the first one and of new uniform estimates. Indeed, as it will turn out at the end of this subsection, after letting $\delta$ to $0$ one recovers estimate \eqref{un_main}.  This estimate can be used to substitute all $\kappa$-dependent estimates used to pass to the limit with $\delta$.

\noindent We start with an observation that according to \eqref{un_delta}, there exist subsequences $\rd,\vu_\delta$ such that
\eq{\label{conv}
&\rd\longrightarrow \rho \qquad \text{weakly-*  in }  L^{\infty}(0,T;L^{K}(\Omega)), \\
& \ud \longrightarrow \vu \qquad \text{weakly in }  L^2(0,T;(W^{1,2}(\Omega))^3).
}
Next, directly from the continuity equation it follows that $\rd$ is equi-continuous with values in $W^{-1,\frac{2K}{K+1}}(\Omega)$. Moreover $\rd$ belongs to $\mathcal{C}([0,T];L^{K}_{\rm weak}(\Omega))$ and is uniformly bounded in $L^\infty(0,T; L^{K}(\Omega))$, thus by the Arzel\`a-Ascoli theorem we verify that
\begin{equation}\label{cvg_dens}
\rd\longrightarrow \rho \qquad \text{in }\ \mathcal{C}([0,T];L^{K}_{\rm weak}(\Omega)).
\end{equation}
Similarly,  we prove that
\begin{equation}\label{cvg_momentum}
\rd\ud \longrightarrow \rho \vu \quad \text{in }\  \mathcal{C}([0,T]; (L_{\rm weak}^{\frac{2K}{K+1}}(\Omega))^3).
\end{equation}
Indeed, from the momentum equation we conclude that $\rd\ud$ is equi-continuous with values in $W^{-1,s}(\Omega)$ with $s=\frac{K+1}{K}$. Moreover,  on account of \eqref{interp} it is uniformly bounded in $L^\infty(0,T;(L^{\frac{2K}{K+1}}(\Omega))^3)$ and belongs to $\mathcal{C}([0,T]; (L_{\rm weak}^{\frac{2K}{K+1}}(\Omega))^3)$, thus the Arzel\`a-Ascoli theorem yields \eqref{cvg_momentum}.

\bigskip

\noindent{\bf Passage to the limit in the continuity and the momentum equations.} 
Convergences \eqref{cvg_dens} and \eqref{cvg_momentum} allow us to pass to the limit in the approximate continuity equation, which is now satisfied in the sense of distributions on $(0,T)\times\Omega$. 

To pass to the limit in the momentum equation we need to justify the convergence in the nonlinear terms. For the convective term we have
\begin{equation}\label{cvg_conv}
\rd\ud\otimes \ud \longrightarrow \rho \vu \otimes \vu \quad \text{weakly in } L^q((0,T)\times \Omega)
\end{equation}
for some $q>1$. It follows from \eqref{interp} which gave us uniform  bound on $\rd |\ud|^2$ in $L^2(0,T;L^p(\Omega))$ 
with $\displaystyle p= \frac{6K}{4K + 3}>1$, the convergence established in \eqref{cvg_momentum} and in
\eqref{conv} together with compact imbeddings.

Finally, thanks to \eqref{equiintegr} we deduce existence of a subsequence such that
\eq{\Pd\to \overline{\pi_{\vep}\lr{\frac{\rho}{\rs}}} \text{ weakly in } L^1((0,T)\times \Omega),}
however, identification of the limit term cannot be done yet, since $\rho_\delta$ is only known to converge weakly to $\rho$. Nevertheless, with these convergences at hand, the passage to the limit in the continuity equation is automatic. Letting  $\delta\rightarrow 0$ in the weak formulation of the momentum equation, we get
\eq{\label{weak_mom_delta}
&\intTO{\rho \vu\cdot\partial_t\varphi }+ \intO{ \vc{m}_0\cdot\varphi(0)} \\
&+\intTO{\rho(\vu\otimes \vu):\na\varphi} -\intTO{\vc{S} : \na\varphi } \\
&+\intTO{\overline{p(\rho)} \Div\varphi} +\intTO{ \overline{\pi_{\vep}\lr{\frac{\rho}{\rs}}}\Div(\rs\vp)}=  0
}
satisfied for any $\varphi\in C^\infty_c\lr{[0,T)\times\Ov{\Omega}}^3$ such that $\varphi\big|_{\bom}=0$.

The final step is to identify 
\eq{\label{identy}
\overline{p(\rho)}=p(\rho), \quad \mbox{and}\quad \overline{\pi_{\kappa}\lr{\frac{\rho}{\rs}}}= \pi_{\kappa}\lr{\frac{\rho}{\rs}}.}
To do that we need to show the strong convergence of the density.

\bigskip

\noindent {\bf Strong convergence of the density.} Proving the strong convergence of the density is a standard difficulty in the study of compressible Navier-Stokes equations (see \cite{PLL}, \cite{NS}). The main ingredients are the theory of renormalized solutions to the continuity equation (whose definition is recalled in Appendix with the main existence result) and the compactness of a quantity called {\it the effective flux}. We follow in this part the standard theory of compressible Navier-Stokes equations.

The renormalized version of the continuity equation \eqref{renormal} with $b(\rd)=\rd\log\rd$ (we extend $\rd$ and $\vu_\delta$ by $0$ outside $\Omega$) reads 
\begin{equation}\label{renorm}
\partial_t(\rd\log\rd) + \Div(\rd\log\rd\ud) = -\rd\Div \ud \quad \text{in}\quad \mathcal{D}'((0,T)\times \mathbb{R}^3).
\end{equation}
We now want to pass to the limit $\delta\to 0$ in the above equality. Due to \eqref{hyp_renorm_2} with $\lambda_1=1$ and the weak-* convergence of $\rd$ in $L^\infty(0,T;L^K(\Omega))$ established in previous paragraph, we have
\begin{eqnarray*}
& b(\rd)\longrightarrow \overline{b(\rho)} \quad \text{weakly-* in }  L^\infty(0,T;L^{\frac{K}{2}}(\Omega)) &\\
& (\rd b'(\rd)-b(\rd))\Div\ud \rightarrow \overline{(\rho b'(\rho)-b(\rho))\Div \vu} \quad \text{weakly in } L^2(0,T;L^{\frac{2K}{K+4}}(\Omega))&
\end{eqnarray*} 
To prove the convergence in the term $\rd\log\rd\,\ud$ we need the following compensated compactness lemma (see for instance \cite{PLL} Lemma 5.1).
\begin{Lemma}\label{lem_compac} Let $g_n$, $h_n$ converge weakly to $g$, $h$ respectively in $L^{p_1}(0,T;L^{p_2}(\Omega))$, $L^{q_1}(0,T;L^{q_2}(\Omega))$ where $1\leq p_1,p_2 \leq \infty$, $\frac{1}{p_1}+ \frac{1}{q_1}=\frac{1}{p_2}+ \frac{1}{q_2}= 1$. We assume in addition that
\begin{itemize}
  \item $\partial_t g_n$ is bounded in $L^1(0,T;W^{-m,1}(\Omega))$ for some $m\geq 0$ independent of $n$.
  \item $\|h_n-h_n(t,\cdot+\xi)\|_{L^{q_1}(L^{q_2})} \rightarrow 0$ as $|\xi|\rightarrow 0$, uniformly in $n$.
\end{itemize}
Then $g_nh_n$ converges to $gh$ in $\mathcal{D}'$.
\end{Lemma}
This result applied to $g_\delta=\rd\log\rd$ and $h_\delta=\ud$  gives the convergence of $\rd\log\rd\ud$ to $\overline{\rho \log\rho}\vu$, and, passing to the limit $\delta \rightarrow 0$ in \eqref{renorm}, we check that
\eq{\partial_t (\overline{\rho \log\rho} )+ \Div(\overline{\rho\log \rho}\vu) =-\overline{\rho \Div \vu}\label{lim_ren}}
in the sense of distributions.

Comparing \eqref{lim_ren} and the renormalized continuity equation satisfied by the limit functions $\rho$, $\vu$ with $b(\rho)=\rho\log\rho$ we obtain
\begin{equation}\label{ineg_log}
\partial_t(\overline{\rho\log\rho}-\rho\log\rho) + \Div[(\overline{\rho\log\rho}-\rho\log\rho)\vu]=-\overline{\rho\Div \vu} +\rho\Div \vu.
\end{equation}
We will exploit this equality to show compactness of the, so called, effective viscous flux.
To derive a key equality for this reasoning, we introduce the the inverse divergence operator
	$\mathcal{A}=\na\lap^{-1}$
and the double Riesz transform 	
	$\mathcal{R}=\na\otimes\na\Delta^{-1}$
specified by \eqref{defA} and \eqref{defR} in Appendix, where we also recall some of their basic properties.
	\bigskip

Applying the operator $\na \Delta^{-1}$ to the approximate continuity equation, we get
\[\partial_t\na\Delta^{-1}(\rd\un_\Omega) + \na \Delta^{-1}\Div(\rd\ud) = 0.\] 
Due to Lemma \ref{LA1}, it follows that $\na\Delta^{-1}(\rd\un_\Omega) \in L^{\alpha+1}(0,T;W^{1,\alpha+1}(\Omega))\cap L^\infty(0,T;W^{1,\alpha}(\Omega))$ and $\partial_t\na\Delta^{-1}(\rd\un_\Omega) \in L^2((0,T)\times\Omega)$.
Thus 
$$\varphi_\delta = \psi\phi\na\Delta^{-1}[\un_\Omega\rd]$$
with $\phi\in \mathcal{D}(\Omega)$ and $\psi\in\mathcal{D}((0,T))$ is an admissible test function for the momentum equation. After straightforward calculation we obtain
\eq{\label{eff_flux_delta}
&\intT{\psi\intO{\phi\rd\left[p(\rd)+\rs\Pd - (2\mu+\lambda)\Div\ud\right]}} \\
& = -\intT{\psi\intO{\left[p(\rd)+\rs\Pd\right] \partial_i\phi\mathcal{A}_i(\rd\un_\Omega)}}  \\
& \quad-\intT{\psi\intO{\phi\Pd \partial_i\rs\mathcal{A}_i(\rd\un_\Omega)}}  \\
& \quad + (\mu+\lambda)\intT{\psi\intO{\partial_i\phi\Div\ud\mathcal{A}_i(\rd\un_\Omega)}}\\
&\quad+ \mu\intT{\psi\intO{\partial_j\phi\partial_ju_\delta^i\mathcal{A}_i(\rd\un_\Omega)}}\\ 
& \quad -\intT{\psi\intO{\partial_j\phi\rd u_\delta^i u_\delta^j\mathcal{A}_i(\rd\un_\Omega)}} \\
&\quad -\intT{\partial_t\psi\intO{\phi\rd u_\delta^i\mathcal{A}_i(\rd\un_\Omega)}} \\
& \quad +\intT{\psi\intO{u_\delta^j\left[\rd\un_\Omega\mathcal{R}_{ij}\big(\rd u_\delta^i\phi\big) -\rd u_\delta^i\phi\mathcal{R}_{ij}(\rd\un_\Omega)\right]}}
}
Using the convergences established above and the properties of operator $\mathcal{A}_i$ (Lemma \ref{LA1}), we may pass to the limit $\delta\to 0$ in all the terms of the above formula except the last one. 
Since $\int_{\mathbb{R}^3}{\mathcal{R}_{ij}(f)g} = \int_{\mathbb{R}^3}{\mathcal{R}_{ij}(g)f}$, for all $ f\in L^r(\mathbb{R}^3),\, g\in L^{r'}(\mathbb{R}^3)$
we deduce that the last term can be rewritten as
\eq{\label{com}
&\intT{\psi\intO{\rd\left[u_\delta^j\mathcal{R}_{ij}\big(\rd u_\delta^i\phi\big) -\mathcal{R}_{ij}\big(\rd u_\delta^i\phi u_\delta^j\big)\right]}} \\
& = \intT{\psi\intO{\rd[u_\delta^j,\mathcal{R}_{ij}](\phi\rd u_\delta^i)}}.
}
Due to \eqref{un_delta} and \eqref{interp}, Lemma \ref{LA3} may be applied to control\begin{equation*}
[\ud,\mathcal{R}](\phi\rd\ud) \text{ in } L^1(0,T;W^{1,q}(\Omega)) \text{ with } \frac{1}{q} = \frac{1}{2}+\frac{\alpha+ 6}{6\alpha}=\frac{2\alpha+3}{3\alpha},
\end{equation*}
so, since $\alpha>3$,  $q\in (1,\frac{3}{2})$.\\
Next, we identify the limit  of $[\ud,\mathcal{R}](\phi\rd\ud)$. Due to \eqref{cvg_conv} we already know that $\mathcal{R}_{ij}(\rd u_\delta^i u_\delta^j)$ converges weakly to $\mathcal{R}_{ij}(\rho u^i u^j)$ thus we only have to justify the weak convergence of $u_\delta^j\mathcal{R}_{ij}(\rd u_\delta^i)$. This easily follows by application of Lemma \ref{lem_compac} to $g_\delta=\mathcal{R}_{ij}(\phi \rd\ud^i)$ and $h_\delta=u_\delta^j$. Thus, Lemma \ref{LA2} yields
\begin{equation}
[u_\delta^j,\mathcal{R}_{ij}](\phi\rd u_\delta^i) \quad\text{converges weakly to }[u^j,\mathcal{R}_{ij}](\phi\rho u^i)\quad \text{in } L^1(0,T,L^p(\Omega)), 
\end{equation}
for $p$ such that $W^{1,q}(\Omega) \hookrightarrow L^p(\Omega)$, namely
\[p < \frac{3q}{3-q} = \frac{3\alpha}{\alpha{+}3}.\]
Observe that since $\alpha>3$ thus $p\leq3/2$ and therefore the integral on the r.h.s. of \eqref{com} is bounded provided $\rho_\delta\in L^\infty(0,T;L^{p'}(\Omega))$ with $p'\geq 3$, which is satisfied for $K\geq3$. Finally, applying once more Lemma \ref{lem_compac} this time with $g_\delta=\rd$ and $h_\delta=[u_\delta^j,\mathcal{R}_{ij}](\phi\rd u_\delta^i) $ we identify the limit of \eqref{com}.
%
Therefore, passing to the limit $\delta\to 0$ in \eqref{eff_flux_delta}, we obtain
\eq{\label{eff_flux_lim}
&\lim_{\delta\to 0}\intT{\psi\intO{\phi\rd\left[p(\rd)+\rs\Pd - (2\mu+\lambda)\Div\ud\right]}} \\
& = -\intT{\psi\intO{\left[\Ov{p(\rho)}+\rs \overline{\pi_{\vep}\lr{\frac{\rho}{\rs}}}\right] \partial_i\phi\mathcal{A}_i(\rd\un_\Omega)}}  \\
& \quad-\intT{\psi\phi\intO{ \overline{\pi_{\vep}\lr{\frac{\rho}{\rs}}}\partial_i\rs\mathcal{A}_i(\rho\un_\Omega)}} \\
& \quad + (\mu+\lambda)\intT{\psi\intO{\partial_i\phi\Div \vu\mathcal{A}_i(\rho\un_\Omega)}}\\
&\quad+ \mu\intT{\psi\intO{\partial_j\phi\partial_j u^i\mathcal{A}_i(\rho\un_\Omega)}}\\
& \quad -\intT{\psi\intO{\partial_j\phi\rho u^i u^j\mathcal{A}_i(\rho\un_\Omega)}} -\intT{\partial_t\psi\intO{\phi\rho u^i\mathcal{A}_i(\rho\un_\Omega)}} \\
& \quad +\intT{\psi\intO{\rho\left[u^j,\mathcal{R}_{ij}\right](\phi\rho u^i)}}.
}
This is to be compared with analogous expression obtained for a limit momentum equation \eqref{weak_mom_delta} with
$\varphi = \psi\phi\na\Delta^{-1}[\un_\Omega\rho]$, where $\phi\in \mathcal{D}(\Omega)$ and $\psi\in\mathcal{D}((0,T))$, we have
\eq{\label{eff_flux2}
&\intT{\psi\intO{\phi\rho\left[\overline{p(\rho)}+ \rs\overline{\pi_{\vep}\lr{\frac{\rho}{\rs}}} - (2\mu+\lambda)\Div \vu\right]}} \\
& = -\intT{\psi\intO{\left[\overline{p(\rho)}+ \rs\overline{\pi_{\vep}\lr{\frac{\rho}{\rs}}}\right] \partial_i\phi\mathcal{A}_i(\rho\un_\Omega)}}  \\
& \quad-\intT{\psi\phi\intO{ \overline{\pi_{\vep}\lr{\frac{\rho}{\rs}}}\partial_i\rs\mathcal{A}_i(\rho\un_\Omega)}} \\
& \quad + (\mu+\lambda)\intT{\psi\intO{\partial_i\phi\Div \vu\mathcal{A}_i(\rho\un_\Omega)}}\\
&\quad+ \mu\intT{\psi\intO{\partial_j\phi\partial_j u^i\mathcal{A}_i(\rho\un_\Omega)}}\\
 & \quad -\intT{\psi\intO{\partial_j\phi\rho u^i u^j\mathcal{A}_i(\rho\un_\Omega)}} -\intT{\partial_t\psi\intO{\phi\rho u^i\mathcal{A}_i(\rho\un_\Omega)}} \\
& \quad +\intT{\psi\intO{\rho\left[u^j,\mathcal{R}_{ij}\right](\phi\rho u^i)}}.
}
Comparing \eqref{eff_flux_lim} with \eqref{eff_flux2} we get
\eq{\label{ef}
&\lim_{\delta\to 0}\intT{\psi\intO{\phi\rd\left[p(\rd)+\rs \Pd- (2\mu+\lambda)\Div\ud\right]}} \\
& \qquad = \intT{\psi\intO{\phi\rho\left[\overline{p(\rho)}
+\rs\overline{\pi_{\kappa}\lr{\frac{\rho}{\rs}}}- (2\mu+\lambda)\Div \vu\right]}}
.}
On one hand, due to monotonicity of $p(\rho)$ we have 
\[\intTO{\phi\big(p(\rd)\rd-\Ov{p(\rho)}\rho\big)} \geq 0 .\]
On the other hand, since $\cdot\mapsto\pi_{\kappa,\delta}(\cdot)$ is non-decreasing for fixed $\delta$ and $\pi_{\kappa,\delta^1}(\cdot) \geq \pi_{\kappa,\delta^2}(\cdot)$ provided $\delta^1\leq \delta^2$, therefore
\eq{ \liminf_{\delta\rightarrow 0} \intTO{\phi\rs\lr{\Pd\rd-\overline{\pi_{\kappa}\lr{\frac{\rho}{\rs}}}\rho}}\\
 \geq 
 \liminf_{\delta\rightarrow 0}\intTO{\phi\rs \lr{\pi_{\kappa,\delta^0}\lr{\frac{\rd}{\rs}}\rd 
 -\overline{\pi_{\kappa}\lr{\frac{\rho}{\rs}}}\rho}}\\
 \geq 
 \liminf_{\delta\rightarrow 0}\intTO{\phi \rs\lr{\Ov{\pi_{\kappa,\delta^0}\lr{\frac{\rho}{\rs}}} -\overline{\pi_{\kappa}\lr{\frac{\rho}{\rs}}}}\rho},
}
where the last inequality is again a consequence of the monotonicity of $\pi_{\delta^0}$ for any fixed $\delta^0$. Now, one can see that the r.h.s. vanishes due to the strong convergence of $\Ov{\pi_{\delta^0}\lr{\frac{\rho}{\rs}}} $ to $\overline{\pi_{\vep}\lr{\frac{\rho}{\rs}}}$ for $\delta^0\to 0$. This follows from the equi-integrability of the singular pressure \eqref{equiintegr}.

These inequalities imply that \eqref{ef} can be reduced to
\[\intTO{\overline{\rho\Div \vu}} \geq \intTO{\rho\Div \vu}.\] 
Coming back to \eqref{ineg_log} and using the convexity of the function $s\mapsto s\log s$ we get
\[\overline{\rho\log \rho} = \rho\log \rho\]
which yields the strong convergence of $\rho_\delta$ in $L^p((0,T)\times\Omega)$, $p<\alpha+1$.\\

Having disposed of the problem of strong convergence of the density, we can identify the limits \eqref{identy} in the momentum equation \eqref{weak_mom_delta}. 

In addition, using Mosco convergence, one can let $\delta\to 0$ also in the energy inequality \eqref{energ_in}, we have

\eq{\label{energ_delta}
\Dt \intOB{\frac{1}{2} \rho |\vu|^2 + \frac{1}{\gamma-1}\rho^\gamma
+ \rho\Gamma_{\kappa}\lr{\frac{\rho}{\rho^*}}} 
+ \intOB{\vc{S}:\nabla\vu} \leq0 
}
satisfied in the sense of distributions on $(0,T)$.\\

\bigskip

\noindent{\bf Upper bound on the limit density.}  To conclude this section, let us prove the uniform bound for limit density $\re$, which is the main advantage of our approximation scheme in comparison with \cite{LiMa}.
Recall that the basic energy estimate \eqref{un_delta} implies in particular that
\[\displaystyle \sup_{t\in[0,T]} \intO{\rho_{\delta} \Gamma_{\kappa,\delta}\lr{\frac{\rd}{\rho^*}}} \leq C.\]
Thus, the upper bound for $\rho_\kappa$ follows from  
\begin{align*}
\int_\Omega{\rd\Gamma_{\kappa,\delta}\left(\frac{\rd}{\rs}\right)\dx} &\geq \int_\Omega{\rd\Gamma_{\kappa,\delta}\left(\frac{\rd}{\rs}\right)\vc{1}_{\{\rd/\rs\geq 1-\delta\}}\dx} \\
& \geq \int_\Omega{\rd\left(\int_0^{1-\delta}{\frac{s^{\alpha-2}}{(1-s)^\beta}\mathrm{d}s}\right)\vc{1}_{\{\rd/\rs\geq 1-\delta\}}\dx} \\
&=\vep\int_\Omega{\rd\left(\int_\delta^{1}{\frac{(1-s)^{\alpha-2}}{s^\beta}\mathrm{d}s}\right)\vc{1}_{\{\rd/\rs\geq 1-\delta\}}\dx} \\
&= \vep\int_\Omega{\rd\left(\sum_{k=0}^{\alpha-2}\begin{pmatrix}\alpha -2 \\ k \end{pmatrix}(-1)^k\int_\delta^1{s^{k-\beta}\mathrm{d}s}\right)\vc{1}_{\{\rd/\rs\geq 1-\delta\}}\dx} \\
& =\frac{\vep}{\beta-1}\intO{\rd\lr{\frac{1}{\delta^{\beta-1}}-1}\vc{1}_{\{\rd/\rs\geq 1-\delta\}}}
\\
&\quad+ \vep\int_\Omega{\rd\left(\sum_{k=1}^{\alpha-2}\begin{pmatrix}\alpha -2 \\ k \end{pmatrix}(-1)^k\int_\delta^1{s^{k-\beta}\mathrm{d}s}\right)\vc{1}_{\{\rd/\rs\geq 1-\delta\}}\dx}
\end{align*}
Similarly to \eqref{ineg_equi} we get
\[\geq \vep\frac{C_1(\vep,\alpha,\beta)}{\delta^{\beta-1}}\left|\left\{\frac{\rd}{\rs}\geq 1-\delta\right\}\right|-C_2(\vep,\alpha,\beta)\]
which implies that $\left|\left\lbrace\dfrac{\rho_\delta}{\rs}\geq 1-\delta\right\rbrace\right| \leq C(\vep)\delta^{\beta-1}$. Recall that $\beta>3$, thus finally, after letting $\delta\to 0$ we obtain $\left|\left\lbrace\dfrac{\rho_\kappa}{\rs}\geq 1\right\rbrace\right|=0$, which implies \eqref{un_main}.
%
%
%

Having obtained this uniform bound, passage to the limit $\kappa\to 0$ is just a repetition of the steps from above and thus,
\begin{equation}\label{rho_bound}
0 \leq \rho_\vep \leq 1,
\end{equation}
 and the first part of Theorem \ref{main} is proved .$\Box$

\section{Recovering of the two-phase system}\label{S_ep}
The purpose of this section is to perform the last limit passage, i.e. $\vep\to 0$ and so to prove the second part of Theorem \ref{main}. 

\bigskip

\noindent {\bf Uniform bounds.} We first summarize what kind of uniform estimates are available at this level of approximation.
Directly from \eqref{energ_delta} it follows that
\eq{
& (\re |\ue|^2)_{\{\vep>0\}}\quad \text{is bounded in } L^{\infty}(0,T;L^1(\Omega)) \\
& (\re )_{\{\vep>0\}} \quad \text{is bounded in } L^{\infty}((0,T)\times \Omega) \\
& \lr{\re \Gamma_\vep\lr{\frac{\re}{\rs}}}_{\{\vep>0\}}\quad  \text{is bounded in } L^{\infty}(0,T;L^1(\Omega)) \\
& (\ue)_{\{\vep>0\}}\quad  \text{ is bounded in } L^2(0,T;(W^{1,2}(\Omega))^3).
}
Next, since the density sequence is uniformly bounded in $L^\infty((0,T)\times \Omega)$, we can test the momentum equation by
\[\phi(t,x) = \frac{\psi(t)}{\rs}\mathcal{B}\left(\re - \Ov{\re}\right),\qquad \psi(t)\in {\mathcal C}^\infty_0((0,T)),\quad \psi\geq0,\]
to get the uniform bound for the pressure
\eq{\label{pressure_bound_eps}
\intTO{\Pn} \leq C,\quad\text{and}\quad {\intTO{\re \Pn}} \leq C.
}

\bigskip

\noindent{\bf Passage to the limit in the continuity and the momentum equations.} Using the approximate  continuity and  momentum equations we may repeat the steps leading to \eqref{cvg_dens} to get
\begin{equation}
\re \longrightarrow \rho \quad\text{ in } \ \mathcal{C}([0,T];L_{\rm weak}^p(\Omega)) \quad \forall 1\leq p <+\infty.
\end{equation}

Moreover, the compensated compactness lemma \ref{lem_compac} allows to justify the convergence of $\re\ue$ and $\re\ue\otimes\ue$ in the sense of distributions to $\rho \vu$ and $\rho \vu\otimes \vu$ respectively, as it was done in \eqref{cvg_momentum} and \eqref{cvg_conv}. These observations are sufficient in order to pass to the limit in the approximate continuity equation in order to obtain \eqref{AC}. 

Further, thanks to \eqref{pressure_bound_eps} we can extract the subsequences such that
\begin{align*}
& \Pn \longrightarrow {\pi} \quad \text{weakly in } \mathcal{M}_+((0,T)\times \Omega), \\
& \re\Pn\longrightarrow {\pi}_1\quad \text{weakly in } \mathcal{M}_+((0,T)\times \Omega).
\end{align*}
This allows us to pass to the limit in the momentum equation \eqref{weak_mom_delta}, we obtain
\eq{\label{weak_mom}
&\intTO{\rho \vu\cdot\partial_t\varphi }+ \intO{ \vc{m}_0\cdot\varphi(0)} \\
&+\intTO{\rho(\vu\otimes \vu):\na\varphi} -\intTO{\vc{S} : \na\varphi } \\
&+\intTO{\overline{p(\rho)} \Div\varphi} +\intTO{ {\pi}\Div(\rs\vp)}=  0,
}
where the last product of the l.h.s. has to be understood in the sense of distribution since $\pi$ is merely a measure. In particular, this justify the regularity imposed on $\rho^*$ (that it belongs too $ \mathcal{C}^1$). \\ 
As in the previous section, we still need to identify the limit $\Ov{p(\rho)}$ but also to prove that the limit measure $\pi$ satisfies the constraint from \eqref{Cons}.

\bigskip

\noindent{\bf Strong convergence of the density.} The strong convergence of the sequence approximating the density will allow us to identify 
\eq{\label{identy2}
\overline{p(\rho)}=p(\rho).}
In order to prove that, we need to derive a variant of effective viscous flux equality \eqref{ef}. The basic idea is the same as previously, we want to test the approximate momentum equation \eqref{pressure_bound_eps} by $\na \Delta^{-1}[\un_\Omega\re]$, which is now bounded in  $L^{\infty}(0,T;W^{1,\infty}(\Omega))$, pass to the limit and compare it with analogous expression obtained for the limit equation \eqref{weak_mom}.
%
Note, however, that we do not have enough regularity on the limit pressure $\pi$ to test the limit momentum equation by  $\na \Delta^{-1}[\un_\Omega\rho]$. 

To justify this step we regularize in time and space the weak limits $\rho$ and $\pi$ by means of standard multipliers. Indeed, taking in \eqref{weak_mom} the supremum over $\vp=\frac{\phi}{\rs}$ for all $\phi \in {\cal D}((0,T)\times\Omega)$ and using the uniform estimates, we can verify that the limit pressure $\pi$ is more regular. 
Indeed, considering the limit momentum equation 
\[\partial_t(\rho \vu) + \Div\,(\rho\vu\otimes \vu) -\Div \,({\bf S}) + \nabla p +\rho^*\nabla \pi = 0 \]
satisfied in the sense of distributions, we can divide this equation by $\rho^*$ and apply the operator $\nabla\Delta^{-1}$ to obtain
\begin{equation}\label{regu_pi}
\pi \in W^{-1,\infty}(0,T;W^{1,2}(\Omega)) + L^p(0,T;L^q(\Omega))\quad p,q>1.\end{equation}
On the other hand, from the continuity equation, we easily get
\[\rho \in \mathcal{C}([0,T];L^p)\cap \mathcal{C}^1([0,T];W^{-1,2}(\Omega))\quad 1\leq p< +\infty.\]
Therefore for $\rho_n = \rho\ast\omega_n$, $\pi_n=\pi\ast \omega_n$, where $\omega_n$ is a mollifying sequence, we have the following convergences
\eq{\label{convol}
&\rho_n\longrightarrow \rho \quad \text{in } \mathcal{C}([0,T]; L^p_{\rm weak}(\Omega)) \cap \mathcal{C}^1([0,T],W^{-1,2}(\Omega)) \\
& \pi_n\longrightarrow \pi \quad \text{in } W^{-1,\infty}(0,T;W^{1,2}(\Omega)) + L^p(0,T; L^q(\Omega)).
}
Then the product $\rho\pi$ can be expressed as
\eq{\label{rpi}
\rho\pi=\rho_n\pi_n+(\rho-\rho_n)\pi_n+\rho(\pi-\pi_n)
}
and we can apply Lemma \ref{lem_compac} to pass to the limit in the r.h.s. 
%
%
With this justification, the same computations as in the previous section, give rise to the effective flux equality
\eqh{
&\intTO{\psi\phi\big(\rs\pi_1 +\overline{\rho p(\rho)} - (\lambda+2\mu)\overline{\rho \Div \vu}\big)} \\
&-  \intTO{\phi\big(\rs\rho\pi + \rho\overline{p(\rho)}- (\lambda+2\mu)\rho \Div \vu\big)}=0,
}
where $\phi \in {\cal D}(\Omega)$, $\psi\in {\cal D}$((0,T)). Therefore
%
\eq{ \label{eff_flux}
&(2\mu+\lambda)\intTO{\psi\phi(\overline{\rho\Div \vu}-\rho \Div \vu)}\\
& = \intTO{\psi\phi\big(\overline{p(\rho)\rho}-\overline{p(\rho)}\rho + \rs\pi_1-\rs\rho\pi\big)} \\
&  \geq \intTO{\psi\phi\rs\big(\pi_1-\rho\pi\big) },
}
where the last inequality is a consequence of monotonicity of $p(\rho)$. 

In order to treat the r.h.s. of \eqref{eff_flux} we show that 
\eq{\pi_1 = \rs\pi\label{c0}} 
a.e. in $(0,T)\times\Omega$. To prove this fact we first observe that
\eq{\label{pi_rho_pi}
 \re \pi_\vep\lr{\dfrac{\re}{\rs}} = & \vep \re \dfrac{\left(\dfrac{\re}{\rho^*}\right)^\alpha}{\left(1-\dfrac{\re}{\rho^*}\right)^\beta} \\
  = & \vep \rho^* \dfrac{\left(\dfrac{\re}{\rho^*}\right)^\alpha}{\left(1-\dfrac{\re}{\rho^*}\right)^\beta} + \vep (\re-\rho^*) \dfrac{\left(\dfrac{\re}{\rho^*}\right)^\alpha}{\left(1-\dfrac{\re}{\rho^*}\right)^\beta} \\
  = & \rs\pi_\vep\lr{\dfrac{\re}{\rs}} -\vep\frac{\rs\lr{\dfrac{\re}{\rs}}^\alpha}{\left(1-{\dfrac{\re}{\rs}}\right)^{\beta-1}}. 
}
Then, letting $\vep\rightarrow 0$ and using the $L^1((0,T)\times\Omega)$ bound on $\pi_\vep$, we show that 
\[\vep\frac{\rs\lr{\dfrac{\re}{\rs}}^\alpha}{\left(1-{\dfrac{\re}{\rs}}\right)^{\beta-1}} \longrightarrow 0 \quad \text{strongly in } L^{\beta/(\beta-1)}((0,T)\times \Omega).\]
Therefore, passing to the limit $\vep\to0$ in \eqref{pi_rho_pi}, we arrive at \eqref{c0}.
Inserting this to \eqref{eff_flux} we may deduce that
\eqh{
(2\mu+\lambda)\intTO{\psi\phi(\overline{\rho\Div \vu}-\rho \Div \vu)}  \geq \intTO{\psi\phi\rs\big(\rs-\rho\big)\pi}  \geq 0,
}
where the last inequality follows from the fact that  $\rho\leq \rs$ a.e. in $(0,T)\times\Omega$.

Note that both pairs $(\rho,\vu)$ and $(\re,\ue)$ satisfy the renormalized continuity equation \eqref{renormal}, thus the above equality implies that taking $b(\re)=\re\log\re$ we pass to the weak limit and compare with the equation on $\rho\log\rho$. The same arguments as in the previous section show that $\re$ converges  to $\rho$ strongly in $L^p((0,T)\times \Omega)$, $1\leq p<+\infty$.

\bigskip

\noindent {\bf Recovery of the congestion constraint.} The strong convergence of the density enables also to identify the limit in the singular pressure, we have
\eq{\pi_1=\rho \pi,\label{c1}}
where the meaning to the product on the r.h.s. is given as in \eqref{rpi}.
Finally, comparing  \eqref{c0} with \eqref{c1} we get that
$$\rs\pi=\pi_1=\rho\pi,$$
which gives the congestion constraint \eqref{Cons}.

\bigskip

\noindent {\bf The divergence free condition.} To justify that the limit triple $(\rho,\vu,\pi)$ is a solution to system \eqref{CM0} in the sense of Definition \eqref{Def1} one has to check that the divergence free condition $\Div\,(\rs\vu)=0$ is satisfied a.e. in $\{\rho=\rs\}$.
We will show that it follows from a certain compatibility between  \eqref{CM0}$_{1}$ and conditions \eqref{CM0}$_{2,4}$ for the limit system. We have  the following generalization of Lemma 2.1 from \cite{LiMa}  to the heterogeneous case. 

\begin{Lemma}\label{LMP}
Let $\vu\in L^2(0,T;(W^{1,2}_0(\Omega))^3)$ and $\rho\in L^2((0,T)\times\Omega)$ such that
$$\pt\rho+\Div(\rho \vu)=0\quad in \ (0,T)\times\Omega,\quad \rho(0)=\rho_0$$
then the following two assertions are equivalent
\begin{itemize}
\item[(i)] $\Div(\rho^* \vu) = 0$ a.e on $\{\rho \geq \rho^*\}$ and $0 \leq \rho_0 \leq \rho^*$,
\item[(ii)] $0\leq\rho(x,t)\leq\rho^*(x)$.
\end{itemize}
\end{Lemma}
\pf We first prove implication $(ii)\to (i)$. Denote $\disp R = \frac{\rho}{\rs}$, then $R$ satisfies the equation
\begin{equation}\label{eq_R}
\partial_t R + \Div(R \vu) + R\vu\cdot\na\log \rs = 0.
\end{equation}
Let $\beta \in \mathcal{C}^1([0,\infty))$, multiplying the previous equation by $\beta'$ we get  
\begin{equation} \label{renorm_R}
\partial_t \beta(R) + \Div(\beta(R) \vu) + (\beta'(R)R -\beta(R)) \Div \vu + \beta'(R)R\vu\cdot\na\log\rs=0.
\end{equation}
Due to assumptions $0 \leq R \leq 1$, using in \eqref{renorm_R} $\beta(R) = R^k$ for any integer $k$,  we obtain
\[\partial_t R^k + \Div(R^k\vu) + (k-1)R^k\Div \vu + kR^k\vu\cdot\na\log\rs = 0.\]
Since $R^k \in L^\infty((0,T)\times\Omega)$, thus $\partial_tR^k\in W^{-1,\infty}((0,T)\times \Omega)$. On the other hand, $|R^k\vu|\leq |R\vu|\in L^\infty(0,T;L^2(\Omega))$, so $\Div(R^k\vu)$ is bounded in $L^\infty(0,T;W^{-1,2}(\Omega))$ and $R^k\Div \vu$ is bounded in $L^2((0,T)\times \Omega)$. As a conclusion, we see that\\
$kR^k(\Div \vu + \vu\cdot\na\log\rs)$ is a bounded distribution. Hence, letting $k$ go to infinity we justify that
\[R^k(\Div \vu + \vu\cdot\na\log\rs)  \rightharpoonup 0 \quad \text{in } \mathcal{D}'((0,T)\times\Omega).\]
However, we also know that
\[R^k(\Div \vu + \vu\cdot\na\log\rs) \to \un_{\{R=1\}}(\Div \vu + \vu\cdot\na\log\rs) \quad \text{a.e. in} \ (0,T)\times\Omega \]
and since $|R^k(\Div \vu + \vu\cdot\na\log\rs)|\leq |\Div \vu + \vu\cdot\na\log\rs|$, where  $\Div \vu,\ \vu\cdot \na\log\rs \in L^2((0,T)\times \Omega)$, we get  
\[ \Div \vu + \vu\cdot\na\log\rs = 0 \quad \text{a.e. on} \quad \{R= 1\}.\]

\medskip 

To prove the implication $(i)\to (ii)$ we set $d(\rho)=\rho-\rho^*$ which satisfies the equation
\eq{\partial_td + \Div(d \vu) + \Div(\rho^* \vu) = 0.
\label{eq_d}}
Let us next consider $b(d)$ the positive part of $d(\rho)$, regularized around $0$ in the following way
\[ b_\eta(s) = \begin{cases} b(s) \quad & \text{if } |s|\geq \eta, \\
                             \frac{1}{4\eta}(s+\eta)^2 \quad & \text{if } |s|< \eta. \end{cases}\]
Multiplying \eqref{eq_d} by $b_\eta'(d)$, we obtain the equation
\eq{\partial_tb_\eta(d) + \Div(b_\eta(d) \vu) + (b_\eta'(d)d-b_\eta(d))\Div \vu + b_\eta'(d)\Div(\rho^* \vu) = 0.\label{eq_ren_d}}
Note that
\[ b_\eta(d)'d-b_\eta(d) = \begin{cases} 0 \quad & \text{if } |\rho-\rho^*|\geq \eta, \\
                             \frac{1}{4\eta}(d^2-\eta^2) \quad & \text{if } |\rho-\rho^*|< \eta, \end{cases}\]
and that for $\eta\to0$, $b_\eta'(d)d-b_\eta(d)$, $b'_\eta(d)$  converge strongly in $L^2((0,T)\times\Omega)$ to 0 and $b'_+(d)$ respectively, where $b_+'$ is defined by                          
\[ b_+'(d) = \begin{cases} & 0 \quad \text{if}\quad \rho< \rho^* \\
                           & 1 \quad \text{if}\quad \rho\geq \rho^* \end{cases}.\]
Thus letting $\eta \rightarrow 0$ in \eqref{eq_ren_d} we get
\begin{equation}\label{eq_renorm_b}
\partial_tb(d) + \Div(b(d) \vu) + b_+'(d)\Div(\rho^* \vu) = 0.
\end{equation}
By assumption, $\Div(\rho^* \vu)=0$ on $\{\rho \geq \rho^*\}$ and $b(d)(t=0)= 0$. Therefore, integrating  \eqref{eq_renorm_b}  over $(0,t)\times\Omega$, we obtain $\intO{b(d)(t)} = \intO{b(d)(0)}= 0$ for $t\geq 0$. Since $b(d)$ is a nonnegative  function, we conclude that $d\leq 0$, which means that $\rho \leq \rs$. 
$\Box$

\bigskip

The above result together with estimate \eqref{rho_bound} guarantee that the second part of Theorem \ref{main} is proven. $\Box$


\section{Appendix}\label{Appendix}
{\bf The Bogovskii operator.} The Bogovskii operator is defined in the following lemma.
\begin{Lemma}[\cite{NS}, lemma 3.17]\label{lem_bog}
Let $\Omega$ be a bounded Lipschitz domain in $\mathbb{R}^3$. Then there exits a linear operator $\mathcal{B}_\Omega=\big(\mathcal{B}_\Omega^1,\mathcal{B}_\Omega^2,\mathcal{B}_\Omega^3)$ with the following properties :
\[\mathcal{B}_\Omega\,: \overline{L^p}(\Omega) \rightarrow (W_0^{1,p}(\Omega))^3, \quad 1<p<\infty;\]
\[\Div\,(\mathcal{B}_\Omega(f))=f \quad \text{a.e. in } \Omega,\,f\in \overline{L^p}(\Omega);\]
\[\|\nabla \mathcal{B}_\Omega(f)\|_{L^p(\Omega)} \leq c(p,\Omega)\|f\|_{L^p(\Omega)}, \quad 1<p<\infty\]
\[\text{If } f=\Div\,(g),\text{ with } g\in L^p(\Omega),\, \Div\,(g)\in L^q(\Omega), \, 1<q<\infty, \text{ then }\]
\[\|\mathcal{B}_\Omega(f)\|_{L^q(\Omega)}\leq c(q,\Omega)\|g\|_{L^q(\Omega)}\]
where $\overline{L^p(\Omega)}=\{f\in L^p(\Omega):\int_\Omega{f}(y){\rm d}y=0\}$.
\end{Lemma}
\label{prop_bog}\begin{prop}[\cite{FN}, Theorem 10.11]
Let $\Omega$ be a bounded Lipschitz domain in $\mathbb{R}^3$, $\mathcal{B}_\Omega$ can be uniquely extended as a bounded linear operator
\[\mathcal{B}_\Omega \,: [\overline{W^{1,p'}}(\Omega)]^* = \{f\in [W^{1,p'}(\Omega)]^*; \langle f, 1 \rangle = 0 \} \rightarrow (L^p(\Omega))^3\]
in such way that
\[-\int_\Omega{\mathcal{B}_\Omega(f)\cdot \nabla v} = \langle f,v\rangle_{\{[\overline{W^{1,p'}}(\Omega)]^*,W^{1,p'}\}}, \quad \forall v \in W^{1,p'}(\Omega);\]
\[\|\mathcal{B}_\Omega(f)\|_{L^p(\Omega)} \leq c\|f\|_{[\overline{W^{1,p'}}(\Omega)]^*}.\]
\end{prop}

\bigskip
\noindent{\bf The Riesz transform.} The inverse divergence operator
	$\mathcal{A}=\na\lap^{-1}$
and the double Riesz transform 	
	$\mathcal{R}=\na\otimes\na\Delta^{-1}$
are defined as
\begin{equation}\label{defA}
		\mathcal{A}_j[v]=\left(\na\lap^{-1}\right)_j v=-\mathcal{F}^{-1}\left(\frac{i\xi_j}{|\xi|^{2}}\mathcal{F}(v)\right),
	\end{equation}
	\begin{equation}\label{defR}
		\mathcal{R}_{i,j}[v]=\partial_i\mathcal{A}_j[v]
		=\left(\na\otimes\na\Delta^{-1}\right)_{i,j} v=\mathcal{F}^{-1}\left(\frac{\xi_i\xi_j}{|\xi|^{2}}\mathcal{F}(v)\right).
	\end{equation}
Here, the inverse Laplacian is identified through the Fourier transform $\mathcal{F}$ and the inverse Fourier transform $\mathcal{F}^{-1}$ as
	\begin{equation*}
		(-\Delta)^{-1}(v)=\mathcal{F}^{-1}\left(\frac{1}{|\xi|^{2}}\mathcal{F}(v)\right).
	\end{equation*}
In what follows we recall some of basic properties of these operators.

	\begin{Lemma}\label{LA1}
	The operator ${\cal R}$ is a continuous linear operator from 
	$L^p(\mathbb{R}^3)$ into $L^p(\mathbb{R}^3)$ 
	for any $1<p<\infty$. In particular, the following estimate holds true:
		$$\|{\cal R}[v]\|_{L^p(\mathbb{R}^3)}\leq c(p)\|v\|_{L^p(\mathbb{R}^3)}\quad for\ all\ v\in L^p(\mathbb{R}^3).$$

	The operator ${\cal A}$ is a continuous linear operator from 
	$L^1(\mathbb{R}^3)\cap L^2(\mathbb{R}^3)$ into $L^2(\mathbb{R}^3)+L^{\infty}(\mathbb{R}^3)$, 
	and from $L^p(\mathbb{R}^3)$ into $L^{\frac{3p}{3-p}}(\mathbb{R}^3)$  for any $1<p<3$. 
	
	Moreover,
		$$\|\na {\cal A}[v]\|_{L^p(\mathbb{R}^3)}\leq C(p)\|v\|_{L^p(\mathbb{R}^3)},\quad 1<p<\infty.$$		
	\end{Lemma} 

The proof of this lemma can be found e.g. in \cite{FN},  Section 10.16. In what follows we present two important properties of commutators involving Riesz operator.
The first result is a straightforward consequence of the {\it Div-Curl} lemma, its proof can be found in \cite{EF2001}, Lemma 5.1.

	\begin{Lemma}\label{LA2}
	Let
		\begin{eqnarray*}
			\vc{V}_{\vep}\rightharpoonup \vc{V}\quad {\rm weakly\ in}\ L^{p}(\mathbb{R}^{3}),\quad 
			r_{\vep}\rightharpoonup r\quad {\rm weakly\ in}\ L^{q}(\mathbb{R}^{3}),
		\end{eqnarray*}
	where
		\begin{equation*}
			\frac{1}{p}+\frac{1}{q}=\frac{1}{s}<1.
		\end{equation*}
	Then
		\begin{equation*}
			\vc{V}_{\vep}\mathcal{R}(r_{\vep})-r_{\vep}\mathcal{R}(\vc{V}_{\vep})
			\rightharpoonup 
			\vc{V}\mathcal{R}(r)-r\mathcal{R}(\vc{V})\quad weakly\ in\ L^{s}(\mathbb{R}^{3}).
		\end{equation*}
	\end{Lemma}

The next lemma can be deduced from the general results of {Baj{\v{s}}anski} and {Coifman} \cite{BC}, and {Coifman} and {Meyer} \cite{CM}.

	\begin{Lemma}\label{LA3}
	Let $w\in W^{1,r}(\mathbb{R}^3)$ and $\vc{V}\in L^p(\mathbb{R}^3)$ be given, where
		$1<r<3$, $1<p<\infty$, $\frac{1}{r}+\frac{1}{p}-\frac{1}{3}<\frac{1}{s}<1.$
	Then for all such $s$ we have
		$$\|{\cal R}[w\vc{V}]-w{\cal R}[\vc{V}]\|_{W^{\alpha,s}(\mathbb{R}^3)}
		\leq c(s,p,r)\|w\|_{W^{1,r}(\mathbb{R}^3)}\|\vc{V}\|_{L^p(\mathbb{R}^3)},$$
	where $\alpha$ is given by
		$\frac{\alpha}{3}=\frac{1}{s}+\frac{1}{3}-\frac{1}{p}-\frac{1}{r}.$
	\end{Lemma}
	
 Here,  $W^{\alpha,s}(\mathbb{R}^3)$ for $\alpha\in(0,\infty) \setminus \mathbb{N}$ denotes the Sobolev-Slobodeckii space (see e.g. \cite{T78}).
The proof can be found in \cite{FN}, Section 10.17.

\bigskip

\noindent{\bf The renormalized continuity equation.} Below we recall the definition of the renormalized solution to the continuity equation.

	\begin{df}\label{df2}
	The pair $(\rho,\vu)$ is called a renormalized solution to the continuity equation, if equation \eqref{CM0}$_1$ holds in $\mathcal{D}'((0,T)\times \mathbb{R}^3)$ provided that $\rho,\vu$ is prolonged by $0$ outside $\Omega$ and equation
		\begin{equation}\label{renormal}
			\pt b(\rho)+\Div \lr{b(\rho)\vu}+\lr{\rho b'(\rho)-b(\rho)}\Div\vu=0,
		\end{equation}
	holds in $\mathcal{D}'((0,T)\times \mathbb{R}^3)$ provided that $\rho,\vu$ is prolonged by $0$ outside $\Omega$, for any function $b\in {\cal C}[0,\infty)\cap {\cal C}^{1}(0,\infty)$, such that $sb'(s)\in {\cal C}[0,\infty)$.
	\end{df}

The following result is a consequence of technique introduced and developed by {DiPerna} and {Lions} \cite{DPL}. 
\begin{Lemma}\label{LX}
Let $\rho \in L^{p}((0,T)\times\Omega)$, $p\geq 2$, $\rho\geq 0$, a.e. in $\Omega$ and $\vu\in L^2(0,T; W^{1,2}(\Omega))$ satisfy  continuity equation 
$$\pt\rho+\Div(\rho\vu)=0$$
 in $\mathcal{D}'((0,T)\times \mathbb{R}^3)$ (prolonged by $0$ outside $\Omega$). Then the renormalized continuity equation \eqref{renormal} holds in $\mathcal{D}'((0,T)\times \mathbb{R}^3)$ for any function  $b$ satisfying
\begin{equation}\label{hyp_renorm}
b\in \mathcal{C}[0,+\infty)\cap \mathcal{C}^1(0,+\infty), \quad |b'(t)|\leq ct^{-\lambda_0} \quad t\in (0,+\infty), \quad \lambda_0<1 
\end{equation}
and growth conditions at infinity
\begin{equation}\label{hyp_renorm_2}
|b'(t)|\leq ct^{\lambda_1}, \quad t\geq 1 \quad \text{where } c>0, \quad -1<\lambda_1\leq \frac{p}{2}-1.
\end{equation}
\end{Lemma}
%
%
%
A general reference here is  \cite{FN}, Section 10.18, see also \cite{NS}.\\

\subsection*{Acknowledgment}
   The authors wish to thank Didier Bresch for suggesting the problem, stimulating
conversations and help during the preparation of the paper.

   The first author acknowledges support from the ANR-13-BS01-0003-01 project DYFICOLTI and
 support from the "Projet Exploratoire 2014 de la cellule Energie du CNRS" Dynafilm.
   The second author was supported by MN grant  IdPlus2011/000661.

\end{document}